# Adaptive Boundary Control of Constant-Parameter Reaction-Diffusion PDEs Using Regulation-Triggered Finite-Time Identification

**Iasson Karafyllis[*], Miroslav Krstic[**] and Katerina Chrysafi[*]**

[*]Dept. of Mathematics, National Technical University of Athens, Zografou Campus, 15780, Athens, Greece,
emails: iasonkar@central.ntua.gr , iasonkaraf@gmail.com,
katyusha5393@gmail.com

[**]Dept. of Mechanical and Aerospace Eng., University of California, San Diego, La Jolla, CA 92093-0411, U.S.A., email: krstic@ucsd.edu

**Abstract**

For parabolic PDEs, we present a new certainty equivalence-based adaptive boundary control scheme with a least-squares identifier of an event-triggering type, where the triggering is based on the size of the regulation error (as opposed to the identifier updates being triggered by the estimation error, or the control changes being triggered by the regulation error). The scheme guarantees exponential convergence of the state to zero in the $L^2$ norm and a finite-time convergence of the parameter estimates to the true values of the unknown parameters. The scheme is developed for a specific benchmark problem with Dirichlet actuation, where the only unknown parameters are the reaction coefficient and the high-frequency gain. For this specific problem, no existing adaptive scheme can handle the unknown high-frequency gain. An illustrative example allows the comparison with other adaptive control design methodologies.

## 1. Introduction

The development of adaptive control for parabolic Partial Differential Equations (PDEs) is an actual necessity due to the wide use of parabolic PDEs in many important control problems. In most cases, some of the parameters appearing in the parabolic PDE are not accurately known. Three different methodologies of adaptive boundary control design for parabolic PDEs are available in the literature (see [19,21,23,32,33,34]):
  i)   Lyapunov-based design,
  ii)  design with passive identifiers, and
  iii) design with swapping identifiers.
Recently, the scope of adaptive controllers has been extended to more complicated cases: parabolic PDEs with input delays (see [11]) and parabolic PDEs with distributed parameters and inputs (see [27]). Moreover, adaptive controllers have been used extensively for hyperbolic PDEs: see [1-10,20].

The purpose of the present work is the development of a novel adaptive boundary control scheme for parabolic PDEs. The proposed methodology is based on the extension to the parabolic infinite-dimensional case of the recently proposed adaptive control scheme in [17] for finite-dimensional systems. It is a certainty-equivalence adaptive scheme with a least-squares, regulation-based identifier. The proposed adaptive boundary control scheme guarantees exponential convergence of the state to zero in the $L^2$ norm. Moreover, the scheme guarantees a finite-time convergence of the parameter estimates to the true values of the unknown parameters. The adaptive scheme is developed for a specific benchmark problem with Dirichlet actuation, where the only unknown parameters are the reaction coefficient and the high-frequency gain. For this specific problem, no other adaptive scheme can handle the unknown high-frequency gain. It should be noticed that the proposed adaptive design can be extended to more complicated cases as well as to systems of parabolic PDEs.

An advantage of the certainty-equivalence adaptive boundary control scheme with regulation-based, least-squares identifier is that it can be combined with all methodologies of static boundary feedback design for parabolic PDEs. More specifically, the proposed scheme can be combined with:
  i) the backstepping design (see [22,31]), and
  ii) the reduced model design (see [12,25]).
The derivation of the adaptive control scheme does not require the knowledge of a Lyapunov functional for the parabolic PDE.

It should be noted that the closed-loop system under the proposed adaptive control scheme is a hybrid infinite-dimensional system. The study of hybrid (event-triggered) distributed-parameter systems has attracted the interest of many researchers during the last decade: see the works [13,30] for parabolic PDEs, the works [14,15,24,29] for hyperbolic PDEs and the works [16,26] for abstract infinite-dimensional systems.

The structure of the present work is as follows: Section 2 describes the adaptive boundary control scheme in detail and explains the intuitive ideas behind the event-trigger and the regulation-based, least-squares identifier. The main result (Theorem 2.2) is stated in Section 2 and its consequences are discussed in detail. Section 3 is devoted to the selection of the nominal feedback. It is shown how the adaptive scheme can be combined with the backstepping design and the reduced model design. An illustrative example is presented in Section 4, where the proposed adaptive scheme is compared with the adaptive controller with the passive identifier as well as with the backstepping controller with known parameter. The proof of the main result is provided in Section 5. The proof of Theorem 2.2 requires the development of certain technical, auxiliary results, which are stated in Section 5 and proved in the Appendix. Finally, the concluding remarks are given in Section 6.

**Notation.** Throughout this paper, we adopt the following notation.
* $\Re_+ := [0, +\infty)$. $Z_+$ denotes the set of all non-negative integers.
* Let $U \subseteq \Re^n$ be a set with non-empty interior and let $\Omega \subseteq \Re$ be a set. By $C^0(U;\Omega)$, we denote the class of continuous mappings on $U$, which take values in $\Omega$. By $C^k(U;\Omega)$, where $k \geq 1$, we denote the class of continuous functions on $U$, which have continuous derivatives of order $k$ on $U$ and take values in $\Omega$.
* For a vector $x \in \Re^m$, we denote by $|x|$ its usual Euclidean norm. For real numbers $\alpha_i$, $i = 1,...,m$, $diag(\alpha_1, \alpha_2,...,\alpha_m)$ denotes the diagonal square matrix with $\alpha_1, \alpha_2,...,\alpha_m$ on its main diagonal. For a matrix $W \in \Re^{m \times m}$, the matrix norm $|W|$ is given by $|W| = \sup\{|Wx| : x \in \Re^m, |x| = 1\}$.
* We use the notation $L^2(0,1)$ for the standard space of the equivalence class of square-integrable, measurable functions defined on $(0,1)$ and $\|f\| = \left(\int_0^1 |f(z)|^2 dz\right)^{1/2} < +\infty$, for $f \in L^2(0,1)$.
* For an interval $I \subseteq \Re_+$, the space $C^0(I; L^2(0,1))$ is the space of continuous mappings $I \ni t \to u[t] \in L^2(0,1)$.
* Let $u : \Re_+ \times [0,1] \to \Re$ be given. We use the notation $u[t]$ to denote the profile of $u$ at certain $t \geq 0$, i.e. $(u[t])(z) = u(t,z)$, for all $z \in [0,1]$.

## 2. The Adaptive Controller

In this section we gradually introduce the adaptive control law. The reader interested in a quick access to the adaptive controller may immediately refer to (2.8), (2.9), (2.10), (2.12), (2.13), (2.19), (2.28) and then resume reading the rest of this section for explanations.

Consider the parabolic equation

$$\frac{\partial u}{\partial t}(t,x) = p\frac{\partial^2 u}{\partial x^2}(t,x) + \theta u(t,x), \quad t > 0, \quad x \in (0,1) \tag{2.1}$$

where $p, c > 0$, $\theta \in \Re$ are constants, with boundary conditions

$$u(t,0) = 0, \text{ for } t > 0 \tag{2.2}$$

$$u(t,1) = cU(t), \text{ for } t > 0 \tag{2.3}$$

where $U(t) \in \Re$ is the control input. The values of the parameters $\theta \in \Re$ (reaction coefficient) and $c > 0$ (high-frequency gain) are unknown and are to be estimated.

Let $\Re \ni \theta \to N(\theta) \in \{1,2,3,...\}$ and $\Re \ni \theta \to (k_1(\theta),...,k_{N(\theta)}(\theta)) \in \Re^{N(\theta)}$ be mappings with the property that for every $\theta \in \Re$, there exist constants $R_\theta, \omega_\theta > 0$ such that for every $u_0 \in L^2(0,1)$, the initial boundary value problem (2.1), (2.2), (2.3) with

$$U(t) = c^{-1}\int_0^1 k(\theta, x)u(t,x)dx \tag{2.4}$$

where

$$k(\theta, x) = \sum_{n=1}^{N(\theta)} k_n(\theta)\phi_n(x), \quad x \in [0,1] \tag{2.5}$$

$$\phi_n(x) = \sqrt{2}\sin(n\pi x), \quad n = 1, 2, ... \tag{2.6}$$

and initial condition $u_0 = u[0]$, has a unique solution $u \in C^0(\Re_+; L^2(0,1)) \bigcap C^1((0,+\infty) \times [0,1])$ with $u[t] \in C^2([0,1])$ for $t > 0$ which also satisfies the following estimate

$$\|u[t]\| \leq R_\theta \exp(-\omega_\theta t)\|u_0\|, \quad t \geq 0 \tag{2.7}$$

We refer to the controller (2.4), (2.5) as the nominal feedback. Notice that the implementation of the nominal feedback requires the knowledge of the values of the parameters $\theta \in \Re$, $c > 0$.

The proposed adaptive scheme can work with any nominal feedback of the form (2.4), (2.5) that guarantees (2.7). There are two ways of designing the required mappings $\Re \ni \theta \to N(\theta) \in \{1,2,3,...\}$ and $\Re \ni \theta \to (k_1(\theta),...,k_{N(\theta)}(\theta)) \in \Re^{N(\theta)}$ with the above properties:

1) the finite-mode approximation of the backstepping design (see [22,31]), and
2) the reduced model design methodology (see [12,25]).

An extensive discussion of the construction of the mappings $\Re \ni \theta \to N(\theta) \in \{1,2,3,...\}$ and $\Re \ni \theta \to (k_1(\theta),...,k_{N(\theta)}(\theta)) \in \Re^{N(\theta)}$ is given in the following section (Section 3).

Our certainty equivalence adaptive controller with regulation-triggered least-squares identification has three different components:
  I) the certainty-equivalence controller,
  II) the event-trigger, and
  III) the least squares identifier.
Each component is described in detail below.

## 2.I. The Certainty-Equivalence Controller

The control action in the interval between two consecutive events is governed by the nominal feedback with the unknown $\theta \in \mathfrak{R}$, $c > 0$ replaced by their estimates $\hat{\theta} \in \mathfrak{R}$, $\hat{c} > 0$ at the beginning of the interval. Moreover, the estimates $\hat{\theta} \in \mathfrak{R}$, $\hat{c} > 0$ of the unknown $\theta \in \mathfrak{R}$, $c > 0$ are kept constant between two consecutive events. In other words, we have

$$U(t) = \left(\hat{c}(\tau_i)\right)^{-1} \int_0^1 k(\hat{\theta}(\tau_i), x) u(t, x) dx, \quad t \in [\tau_i, \tau_{i+1}), i \in Z_+ \tag{2.8}$$

$$\hat{\theta}(t) = \hat{\theta}(\tau_i), \quad \hat{c}(t) = \hat{c}(\tau_i), \quad t \in [\tau_i, \tau_{i+1}), i \in Z_+ \tag{2.9}$$

where $\{\tau_i \geq 0\}_{i=0}^{\infty}$ is the sequence of times of the events that satisfies

$$\begin{aligned} \tau_{i+1} &= \min\left(\tau_i + T, r_i\right), \quad i \in Z_+ \\ \tau_0 &= 0 \end{aligned} \tag{2.10}$$

where $T > 0$ is a positive constant (one of the tunable parameters of the proposed scheme) and $r_i > \tau_i$ is a time instant determined by the event trigger.

## 2.II. The Event-Trigger

The proposed event trigger is based on the evolution of the regulation of the state.

Let $a > 0$ be a positive constant (one of the tunable parameters of the proposed scheme). The event trigger sets $r_i > \tau_i$ to be the smallest value of time $t > \tau_i$ for which

$$\|u[t]\| = R_{\hat{\theta}(\tau_i)}(1+a)\|u[\tau_i]\| \tag{2.11}$$

where $u[t]$ denotes the solution of (2.1), (2.2), (2.3) with (2.8) and $R_{\hat{\theta}} > 0$ is the coefficient involved in (2.7). For the case that a time $t > \tau_i$ satisfying (2.11) does not exist, we set $r_i = +\infty$. For the case $u[\tau_i] = 0$ we set $r_i := \tau_i + T$.

Formally, the event trigger is described by the equations:

$$r_i := \inf\left\{t > \tau_i : \|u[t]\| = R_{\hat{\theta}(\tau_i)}(1+a)\|u[\tau_i]\|\right\}, \text{ for } u[\tau_i] \neq 0 \tag{2.12}$$

$$r_i := \tau_i + T, \text{ for } u[\tau_i] = 0 \tag{2.13}$$

## 2.III. The Least-Squares Identifier

The description of the regulation-triggered adaptive control scheme is completed by the parameter update law, which is activated at the times of the events.

Assuming a sufficiently regular solution, we notice that by virtue of (2.1), (2.2), (2.3), we get for all $\tau > 0$ and $n = 1, 2, \ldots$:

$$\frac{d}{d\tau} \int_0^1 \sin(n\pi x) u(\tau, x) dx = -(-1)^n p n \pi c U(\tau) + (\theta - n^2 \pi^2 p) \int_0^1 \sin(n\pi x) u(\tau, x) dx, \tag{2.14}$$

By virtue of (2.14), for every $t, s \geq 0$ and $n = 1, 2, \ldots$, the following equation holds:

$$f_n(t, s) = \theta g_n(t, s) + c j_n(t, s) \tag{2.15}$$

where

$$f_n(t, s) := \int_0^1 \sin(n\pi x)\left(u(t, x) - u(s, x)\right) dx + n^2 \pi^2 p \int_s^t \int_0^1 \sin(n\pi x) u(\tau, x) dx d\tau \tag{2.16}$$

$$g_n(t,s) := \int_s^t \int_0^1 \sin(n\pi x) u(\tau, x) dx d\tau, \quad j_n(t,s) := -(-1)^n pn\pi \int_s^t U(\tau) d\tau \tag{2.17}$$

Let $\tilde{N} \geq 1$ be an (arbitrary; the last of the tunable parameters of the proposed scheme) positive integer. Define for every $i \in Z_+$ and $n = 1, 2, ...$ the function $h_{i,n} : \Re^2 \to \Re_+$ by the formula

$$h_{i,n}(\vartheta_1, \vartheta_2) := \int_{\mu_{i+1}}^{\tau_{i+1}} \int_{\mu_{i+1}}^{\tau_{i+1}} \left( f_n(t,s) - \vartheta_1 g_n(t,s) - \vartheta_2 j_n(t,s) \right)^2 dsdt = 0 \tag{2.18}$$

where

$$\mu_{i+1} := \min \left\{ \tau_j : j \in \{0,...,i\}, \tau_j \geq \tau_{i+1} - \tilde{N}T \right\}. \tag{2.19}$$

The time $\mu_{i+1}$ defined in (2.19) is the time of a past event and is going to play a significant role in what follows: the estimations $\hat{\theta}(\tau_{i+1})$, $\hat{c}(\tau_{i+1})$ at $t = \tau_{i+1}$ are performed by using the measurements of the plant state on the interval $[\mu_{i+1}, \tau_{i+1}]$. Definition (2.19) guarantees that the intervals $[\mu_{i+1}, \tau_{i+1}]$ for $i \in Z_+$ will eventually have large enough length.

It follows from (2.15), (2.18) that for every $i \in Z_+$ and $n = 1, 2, ...$ the function $h_{i,n}(\vartheta_1, \vartheta_2)$ has a global minimum at $(\vartheta_1, \vartheta_2) = (\theta, c)$ with $h_{i,n}(\theta, c) = 0$. Consequently, we get from Fermat's theorem that the following equations hold for every $i \in Z_+$ and $n = 1, 2, ...$:

$$H_{n,1}(\tau_{i+1}, \mu_{i+1}) = \theta Q_{n,1}(\tau_{i+1}, \mu_{i+1}) + c Q_{n,2}(\tau_{i+1}, \mu_{i+1}), \tag{2.20}$$

$$H_{n,2}(\tau_{i+1}, \mu_{i+1}) = \theta Q_{n,2}(\tau_{i+1}, \mu_{i+1}) + c Q_{n,3}(\tau_{i+1}, \mu_{i+1}) \tag{2.21}$$

where

$$H_{n,1}(\tau_{i+1}, \mu_{i+1}) := \int_{\mu_{i+1}}^{\tau_{i+1}} \int_{\mu_{i+1}}^{\tau_{i+1}} f_n(t,s) g_n(t,s) dsdt, \tag{2.22}$$

$$H_{n,2}(\tau_{i+1}, \mu_{i+1}) := \int_{\mu_{i+1}}^{\tau_{i+1}} \int_{\mu_{i+1}}^{\tau_{i+1}} f_n(t,s) j_n(t,s) dsdt \tag{2.23}$$

$$Q_{n,1}(\tau_{i+1}, \mu_{i+1}) := \int_{\mu_{i+1}}^{\tau_{i+1}} \int_{\mu_{i+1}}^{\tau_{i+1}} g_n^2(t,s) dsdt, \tag{2.24}$$

$$Q_{n,2}(\tau_{i+1}, \mu_{i+1}) := \int_{\mu_{i+1}}^{\tau_{i+1}} \int_{\mu_{i+1}}^{\tau_{i+1}} g_n(t,s) j_n(t,s) dsdt, \tag{2.25}$$

$$Q_{n,3}(\tau_{i+1}, \mu_{i+1}) := \int_{\mu_{i+1}}^{\tau_{i+1}} \int_{\mu_{i+1}}^{\tau_{i+1}} j_n^2(t,s) dsdt. \tag{2.26}$$

Define the following set in the parameter space:

$$S_i := \left\{ (\vartheta_1, \vartheta_2) \in \Re^2 : \begin{array}{l} H_{n,1}(\tau_{i+1}, \mu_{i+1}) = \vartheta_1 Q_{n,1}(\tau_{i+1}, \mu_{i+1}) + \vartheta_2 Q_{n,2}(\tau_{i+1}, \mu_{i+1}) \\ H_{n,2}(\tau_{i+1}, \mu_{i+1}) = \vartheta_1 Q_{n,2}(\tau_{i+1}, \mu_{i+1}) + \vartheta_2 Q_{n,3}(\tau_{i+1}, \mu_{i+1}) \end{array}, n = 1, 2, ... \right\} \tag{2.27}$$

Equations (2.20) and (2.21) imply that $(\theta, c) \in S_i$. If $S_i$ is a singleton then a projection on $S_i$ is nothing else but the least-squares estimate of the unknown vector of parameters $(\theta, c)$ on the interval $[\mu_{i+1}, \tau_{i+1}]$.

The following lemma clarifies the form of the set $S_i$ in the parameter space. Its proof can be found in the Appendix.

**Lemma 2.1:** Let $\tau_{i+1} > \mu_{i+1} \geq 0$, $p, c > 0$, $\theta \in \Re$ be given constants and let $u \in C^0\left([\mu_{i+1}, \tau_{i+1}]; L^2(0,1)\right)$, $U \in C^0\left((\mu_{i+1}, \tau_{i+1})\right)$ be given mappings that satisfy (2.15), (2.16), (2.17) for every $t, s \in [\mu_{i+1}, \tau_{i+1}]$ and $n = 1, 2, ...$. Then the following implications hold:
  i) If $Q_{n,1}(\tau_{i+1}, \mu_{i+1}) = 0$ for $n = 1, 2,...$ then $u \equiv 0$ and $S_i = \Re^2$.
  ii) If $S_i = \Re^2$ then $u \equiv 0$ and $Q_{n,1}(\tau_{i+1}, \mu_{i+1}) = 0$ for $n = 1, 2,...$.
  iii) If $S_i \neq \{(\theta, c)\}$ and $S_i \neq \Re^2$ then $U(t) \equiv 0$ and $S_i := \left\{(\theta, \vartheta_2) \in \Re^2 : \vartheta_2 \in \Re \right\}$.

Lemma 2.1 implies that the projection $proj_{S_i}\left((\hat{\theta}(\tau_i), \hat{c}(\tau_i))\right)$ of $(\hat{\theta}(\tau_i), \hat{c}(\tau_i)) \in \Re^2$ on the set $S_i$ can take three values: (i) $(\hat{\theta}(\tau_i), \hat{c}(\tau_i))$ (when $S_i = \Re^2$), (ii) $(\theta, \hat{c}(\tau_i))$ (when $S_i \neq \{(\theta, c)\}$ and $S_i \neq \Re^2$), and (iii) $(\theta, c)$ (when $S_i = \{(\theta, c)\}$). We can therefore define the following parameter update law

$$(\hat{\theta}(\tau_{i+1}), \hat{c}(\tau_{i+1})) = proj_{S_i}\left((\hat{\theta}(\tau_i), \hat{c}(\tau_i))\right) \tag{2.28}$$

Notice that if $\hat{c}(\tau_i) > 0$ then $\hat{c}(\tau_{i+1}) > 0$.

## 2.IV. Main Result: System Properties Under Adaptive Control

Our main result guarantees global exponential regulation of the state $u$ to zero in the $L^2(0,1)$ norm and is stated next.

**Theorem 2.2:** Let $T, a > 0$ be positive constants and let $\tilde{N} \geq 1$ be a positive integer. Then there exists a family of constants $M_{\theta, \hat{\theta}, c, \hat{c}} > 0$, $\omega_\theta > 0$ parameterized by $\theta \in \Re$, $c > 0$, $\hat{\theta} \in \Re$, $\hat{c}_0 > 0$, such that for every $\theta \in \Re$, $c > 0$, $u_0 \in L^2(0,1)$, $\hat{\theta}_0 \in \Re$, $\hat{c}_0 > 0$ the initial-boundary value problem (2.1), (2.2), (2.3) with (2.8), (2.9), (2.10), (2.12), (2.13), (2.19), (2.28) and initial conditions $u[0] = u_0$, $\hat{\theta}(0) = \hat{\theta}_0$, $\hat{c}(0) = \hat{c}_0$ has a unique solution, in the sense that there exist unique mappings $\hat{\theta}, \hat{c} : \Re_+ \to \Re$, $u \in C^0\left(\Re_+; L^2(0,1)\right)$ satisfying $u \in C^1(I \times [0,1])$, $u[t] \in C^2([0,1])$ for $t > 0$, where $I = \Re_+ \setminus \{\tau_i \geq 0, i = 0,1,2,...\}$, which also satisfy (2.8), (2.9), (2.10), (2.12), (2.13), (2.19), (2.28), $u[0] = u_0$, $\hat{\theta}(0) = \hat{\theta}_0$, $\hat{c}(0) = \hat{c}_0$ and

$$\frac{\partial u}{\partial t}(t,x) = p\frac{\partial^2 u}{\partial x^2}(t,x) + \theta u(t,x), \text{ for all } (t,x) \in I \times (0,1) \tag{2.29}$$

$$u(t,0) = u(t,1) - cU(t) = 0, \text{ for all } t \in I \tag{2.30}$$

Moreover, the estimate $\|u[t]\| \leq M_{\theta, \hat{\theta}_0, c, \hat{c}_0} \exp(-\omega_\theta t)\|u_0\|$ holds for all $t \geq 0$. Finally, if $u_0 \neq 0$ then $\hat{\theta}(t) = \theta$ for all $t \geq \tau_1$ and if there exists $t \geq 0$ with $U(t) \neq 0$ then $\hat{c}(t) = c$ for all $t \geq \tau_2$.

It is clear that the proposed regulation-triggered adaptive scheme guarantees global exponential convergence to $0 \in L^2(0,1)$ with exactly the same convergence rate as the nominal feedback. Moreover, the proposed regulation-triggered adaptive scheme guarantees finite-time convergence of the parameter estimates to the exact actual values of the parameters (under certain conditions).

A possible question that may arise at this point is whether the value of the integer $N(\hat{\theta}(t))$ is bounded or not. The proof of Theorem 2.2 shows that $N(\hat{\theta}(t))$ takes only two values when $u_0 \neq 0$: the value $N(\hat{\theta}_0)$ for $t \in [0, \tau_1)$ and the value $N(\theta)$ for $t \geq \tau_1$ that corresponds to the actual value of $\theta$. This is a direct consequence of the fact that $\hat{\theta}(t) = \hat{\theta}_0$ for $t \in [0, \tau_1)$ and $\hat{\theta}(t) = \theta$ for $t \geq \tau_1$ when $u_0 \neq 0$, i.e., the parameter estimation error for $\theta$ becomes zero after the time of the first event (single-trigger convergence). Another possible question that may arise is whether Zeno behavior is possible (i.e., whether the sequence of the times of the events has a finite limit). The proof of Theorem 2.2 shows that Zeno behavior is not possible. More specifically, we have $\tau_i = \tau_2 + (i-2)T$ for all $i \geq 2$ when $u_0 \neq 0$.

## 3. Construction of Nominal Feedback

A nominal feedback law of the form (2.4), (2.5) may be constructed by using two different methodologies.

**1)** The finite-mode approximation of the backstepping design: Theorem 2 in [31] guarantees for every $\theta > 0$, $\beta \geq 0$ the existence of functions $K_\theta, L_\theta \in C^2\left([0,1]^2\right)$ such that the Volterra transformation

$$v(t,z) = u(t,z) - \int_0^z K_\theta(z,s)u(t,s)ds, \text{ for all } (t,z) \in \Re_+ \times [0,1] \quad (3.1)$$

with inverse

$$u(t,z) = v(t,z) + \int_0^z L_\theta(z,s)v(t,s)ds, \text{ for all } (t,z) \in \Re_+ \times [0,1] \quad (3.2)$$

where

$$K_\theta(z,s) := -\frac{\theta+\beta}{p} s \frac{I_1\left(\sqrt{\frac{\theta+\beta}{p}(z^2-s^2)}\right)}{\sqrt{\frac{\theta+\beta}{p}(z^2-s^2)}}, \quad L_\theta(z,s) := -\frac{\theta+\beta}{p} s \frac{J_1\left(\sqrt{\frac{\theta+\beta}{p}(z^2-s^2)}\right)}{\sqrt{\frac{\theta+\beta}{p}(z^2-s^2)}} \text{ for } 0 \leq s < z \leq 1 \quad (3.3)$$

maps the solutions of (2.1), (2.2), (2.3) to the solutions of

$$\frac{\partial v}{\partial t}(t,z) - p\frac{\partial^2 v}{\partial z^2}(t,z) + \beta v(t,z) = 0, \quad (3.4)$$

$$v(t,0) = v(t,1) - U(t) + \int_0^1 \tilde{k}_\theta(s)u(t,s)ds = 0, \quad (3.5)$$

where $\tilde{k}_\theta(z) = K_\theta(1,z)$ for all $z \in [0,1]$. We also define $\tilde{k}_\theta(z) \equiv 0$ for $\theta \leq 0$. Thus for every $\theta \in \Re$, $\beta \geq 0$, we can guarantee that the closed-loop system (2.1), (2.2), (2.3), with $U(t) = \int_0^1 \tilde{k}_\theta(s)u(t,s)ds$ is exponentially stable in the $L^2$ norm.

Following the analysis in the Appendix of [18], we are in a position to guarantee that for every $\theta > 0$, $\sigma \in \left(0, \pi^2 + \frac{\beta}{p}\right)$, $B > 0$ the solution of the initial-boundary value problem (2.1), (2.2), (2.3) with $u[0] = u_0 \in L^2(0,1)$ satisfies the estimate:

$$\|v[t]\| \leq G\exp(-\sigma p t)\|v_0\| + \gamma \sup_{0 \leq s \leq t}\left(\left|U(s) - \int_0^1 \tilde{k}_\theta(x)u(s,x)dx\right|\exp(-\sigma p(t-s))\right), \text{ for } t \geq 0 \quad (3.6)$$

where $v[t]$ is given by (3.1), $G := \sqrt{1+B^{-1}}$ and

$$\gamma := \sqrt{1+B}\begin{cases} \dfrac{(\pi^2+\mu^2)\sqrt{\sinh(2\mu)-2\mu}}{2\sqrt{\mu}(\pi^2+\mu^2-\sigma)\sinh(\mu)} & \text{if } \mu > 0 \\ \dfrac{\pi^2}{\sqrt{3}(\pi^2-\sigma)} & \text{if } \mu = 0 \end{cases} \quad (3.7)$$

where $\mu := \sqrt{\beta/p}$. In this case we may select $N(\theta) \geq 1$ to be an integer sufficiently large so that $\gamma \tilde{L}_\theta \|\tilde{k}_\theta - h_\theta\| < 1$, where $\tilde{L}_\theta := 1 + \left(\int_0^1\left(\int_0^z |L_\theta(z,s)|^2 ds\right)dz\right)^{1/2}$ and $h_\theta(z) := \sum_{n=1}^{N(\theta)} k_n(\theta)\phi_n(z)$ for $z \in [0,1]$ with

$$k_n(\theta) = \int_0^1 \tilde{k}_\theta(x)\phi_n(x)dx, \text{ for } n = 1,...,N(\theta). \quad (3.8)$$

Indeed, using (3.6) and inequality $\gamma \tilde{L}_\theta \|\tilde{k}_\theta - h_\theta\| < 1$ in conjunction with transformations (3.1), (3.2), we are in a position to show that the solution of the closed-loop system (2.1), (2.2), (2.3), with (2.4), (2.5) satisfies estimate (2.7) with $R_\theta := \dfrac{G\tilde{L}_\theta \tilde{K}_\theta}{1 - \gamma \|\tilde{k}_\theta - h_\theta\| \tilde{L}_\theta}$, $\omega_\theta := \sigma p$ and $\tilde{K}_\theta = 1 + \left( \int_0^1 \left( \int_0^z |K_\theta(z,s)|^2 \, ds \right) dz \right)^{1/2}$.

Moreover, using definitions (3.3), we obtain the inequalities $\tilde{L}_\theta \leq 1 + \dfrac{\theta + \beta}{4p\sqrt{3}}$, $\tilde{K}_\theta \leq 1 + \dfrac{1}{2}\sqrt{\dfrac{\theta + \beta}{3p}} I_1\left(\sqrt{\dfrac{\theta + \beta}{p}}\right)$. The previous inequalities in conjunction with (3.7) and the fact that $\|\tilde{k}_\theta - h_\theta\| = \sqrt{\int_0^1 \tilde{k}_\theta^2(x)dx - \sum_{n=1}^{N(\theta)} k_n^2(\theta)}$ (a consequence of (3.8) and definition $h_\theta(z) := \sum_{n=1}^{N(\theta)} k_n(\theta)\phi_n(z)$) allow an estimate of the magnitude of $R_\theta$ and an estimation of how large the integer $N(\theta) \geq 1$ must be selected. This is shown in Table 1 for the case $\beta = 0$, $p = \sigma = 1$, $B = 1/10$.

**Table 1:** Values of the integer $N(\theta) \geq 1$, which satisfy the inequality $\gamma \tilde{L}_\theta \|\tilde{k}_\theta - h_\theta\| < 1$ with $\beta = 0$, $p = \sigma = 1$, $B = 1/10$.

| $\theta$ | $N(\theta)$ |
|---|---|
| 0.1 | 1 |
| 3 | 1 |
| 5 | 2 |
| 6 | 3 |
| 7 | 5 |
| 8 | 7 |
| 9 | 10 |
| $\pi^2$ | 13 |
| 10 | 14 |
| 11 | 19 |
| 12 | 25 |

For $\theta \leq 0$ (recall that $\tilde{k}_\theta(z) \equiv 0$), we define $N(\theta) = 1$ and we still require (3.8) to hold. In this way, we are in a position to show that the solution of the closed-loop system (2.1), (2.2), (2.3), with $\theta \leq 0$, (2.4), (2.5) satisfies estimate (2.7) with $R_\theta := 1$, $\omega_\theta := \pi^2 p$.

**2)** The reduced model design methodology (see [12,25]): The integer $N(\theta) \geq 1$ is selected to be sufficiently large so that $p(N(\theta)+1)^2 \pi^2 > \theta$. Define

$$g(\theta) = \begin{bmatrix} g_1 \\ \vdots \\ g_{N(\theta)} \end{bmatrix} \quad (3.9)$$

where

$$g_n := -p\sqrt{2}(-1)^n n\pi \neq 0 \quad (n = 1, \ldots, N(\theta)) \quad (3.10)$$

The pair of matrices $-p\pi^2 diag(1,\ldots,N^2(\theta)) + \theta I_{N(\theta)}$ and $g(\theta)$ is controllable (see [18]). The mapping $\Re \ni \theta \to (k_1(\theta), \ldots, k_{N(\theta)}(\theta)) \in \Re^{N(\theta)}$ is constructed so that the matrix $-p\pi^2 diag(1,\ldots,N^2(\theta)) + \theta I_{N(\theta)} + g(\theta)\left[k_1(\theta),\ldots,k_{N(\theta)}(\theta)\right]$ is Hurwitz for every $\theta \in \Re$.

## 4. Illustrative Example

In this section we consider the control system (2.1), (2.2), (2.3) with $p = c = 1$ and $\theta = 11$. For these parameter values, the open-loop system (2.1), (2.2), (2.3) with $U(t) \equiv 0$ is unstable with exponentially growing solutions. We assume that the value of the high-frequency gain $c = 1$ is known and the only unknown parameter is the reaction coefficient $\theta$. In this case (known $c > 0$), instead of (2.8) we use the equation

$$U(t) = c^{-1} \int_0^1 k(\hat{\theta}(\tau_i), x) u(t, x) dx, \quad t \in [\tau_i, \tau_{i+1}), i \in Z_+ \tag{4.1}$$

and instead of (2.28) we use the equations

$$\hat{\theta}(\tau_{i+1}) = \begin{cases} \hat{\theta}(\tau_i) & \text{if } Q_n(\tau_{i+1}, \mu_{i+1}) = 0 \text{ for } n = 1, 2, \ldots \\ \dfrac{H_{m(\tau_i)}(\tau_{i+1}, \mu_{i+1})}{Q_{m(\tau_i)}(\tau_{i+1}, \mu_{i+1})} & \text{if otherwise} \end{cases} \tag{4.2}$$

where

$$m(\tau_{i+1}) := \min \{ n \geq 1 : Q_n(\tau_{i+1}, \mu_{i+1}) > 0 \}. \tag{4.3}$$

$$H_n(\tau_{i+1}, \mu_{i+1}) := \int_{\mu_{i+1}}^{\tau_{i+1}} \int_{\mu_{i+1}}^{\tau_{i+1}} f_n(t, s) g_n(t, s) ds dt \tag{4.4}$$

$$Q_n(\tau_{i+1}, \mu_{i+1}) := \int_{\mu_{i+1}}^{\tau_{i+1}} \int_{\mu_{i+1}}^{\tau_{i+1}} g_n^2(t, s) ds dt \tag{4.5}$$

$$\begin{aligned} f_n(t, s) &:= \int_0^1 \sin(n\pi x) \big( u(t, x) - u(s, x) \big) dx \\ &+ (-1)^n cpn\pi \int_s^t U(\tau) d\tau + n^2 \pi^2 p \int_s^t \int_0^1 \sin(n\pi x) u(\tau, x) dx d\tau \end{aligned} \tag{4.6}$$

and $g_n(t, s) := \int_s^t \int_0^1 \sin(n\pi x) u(\tau, x) dx d\tau$ is given by (2.17).

We focus on the application of the backstepping design (explained in the previous section) with $\beta = 0$. The simulations were made by using a finite-difference scheme with 100 spatial grid points. All integrals were evaluated using the trapezoid rule.

Theorem 2 in [31] guarantees that the feedback law

$$U(t) = -11 \int_0^1 s \frac{I_1\left(\sqrt{11(1 - s^2)}\right)}{\sqrt{11(1 - s^2)}} u(t, s) ds \tag{4.7}$$

achieves global exponential stability of the equilibrium point $0 \in L^2(0,1)$ in the $L^2$ norm for the closed-loop system (2.1), (2.2), (2.3) with (4.7). The analysis of the previous section showed that the feedback law (2.4), (2.5) with $N = 20$

$$k_n = -11\sqrt{2} \int_0^1 s \frac{I_1\left(\sqrt{11(1 - s^2)}\right)}{\sqrt{11(1 - s^2)}} \sin(n\pi s) ds, \text{ for } n = 1, \ldots, 20. \tag{4.8}$$

achieves global exponential stability of the equilibrium point $0 \in L^2(0,1)$ in the $L^2$ norm for the closed-loop system (2.1), (2.2), (2.3) with (2.4), (2.5).

The difference between the closed-loop system (2.1), (2.2), (2.3) with (4.7) and the closed-loop system (2.1), (2.2), (2.3) with (2.4), (2.5) was tested numerically. A slight difference appears in the initial transient period, but after this initial transient period there is no visible difference in the response. This is shown in Fig. 1 for the initial condition $u_0(x) = \sqrt{2}\sin(\pi x) + x^2 - x^3$, $x \in [0,1]$.

Next, we develop our regulation-triggered adaptive controller. Following the estimates provided in the previous section, the regulation-triggered adaptive controller is given by (2.9), (2.10), (2.13), (2.19), (2.20), (4.2), (4.3) and

$$U(t) = \sqrt{2} \sum_{n=1}^{N(\hat{\theta}(\tau_i))} k_n(\hat{\theta}(\tau_i)) \int_0^1 \sin(n\pi x) u(t,x) dx, \quad t \in [\tau_i, \tau_{i+1}), i \in \mathbb{Z}_+ \tag{4.9}$$

$$k_n(\vartheta) = -\vartheta\sqrt{2} \int_0^1 s \frac{I_1\left(\sqrt{\vartheta(1-s^2)}\right)}{\sqrt{\vartheta(1-s^2)}} \sin(n\pi s) ds, \text{ for } \vartheta > 0, \ n = 1, 2, \dots. \tag{4.10}$$

$$r_i := \inf\left\{ t > \tau_i : \|u[t]\| = (1+a)M\sqrt{11}\left(4\sqrt{3} + \hat{\theta}(\tau_i)\right)(\pi^2 - 1)\left(\sqrt{3} + \frac{\sqrt{\hat{\theta}(\tau_i)}}{2} I_1\left(\sqrt{\hat{\theta}(\tau_i)}\right)\right)\|u[\tau_i]\| \right\},$$

for $u[\tau_i] \neq 0$ \hfill (4.11)

where $M > 1$, $a > 0$, $\tilde{N} \geq 1$, $T > 0$ are tunable parameters and $N(\vartheta) \geq 1$ is defined for all $\vartheta > 0$ to be the smallest integer for which $12(\pi^2 - 1) - \sqrt{\frac{11}{10}}\pi^2\left(4\sqrt{3} + \vartheta\right)\sqrt{\int_0^1 \tilde{k}_\vartheta^2(x)dx - \sum_{n=1}^{N(\vartheta)} k_n^2(\vartheta)} \geq \frac{1}{M}$ with

$\tilde{k}_\vartheta(s) := -\vartheta s \dfrac{I_1\left(\sqrt{\vartheta(1-s^2)}\right)}{\sqrt{\vartheta(1-s^2)}}$ for $s \in [0,1)$. The quantity appearing in the right hand side of (4.11) is an upper bound of the quantity $\dfrac{G\tilde{L}_\theta \tilde{K}_\theta}{1 - \gamma\|\tilde{k}_\theta - h_\theta\|\tilde{L}_\theta}$, obtained by using (3.7), the inequalities $\tilde{L}_\theta \leq 1 + \dfrac{\theta}{4\sqrt{3}}$,

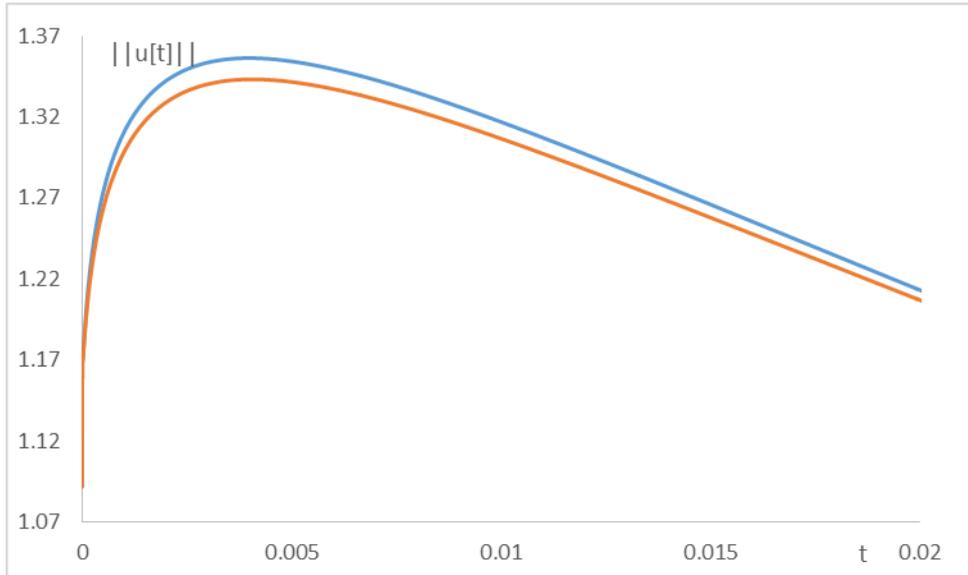

**Fig. 1:** The red curve shows the evolution of $\|u[t]\|$ for the closed-loop system (2.1), (2.2), (2.3) with (4.7). The blue curve shows the evolution of $\|u[t]\|$ for the closed-loop system (2.1), (2.2), (2.3) with (2.4), (2.5) and (4.8).

$\tilde{K}_\theta \leq 1 + \frac{1}{2}\sqrt{\frac{\theta}{3}} I_1(\sqrt{\theta})$, the fact that $\|\tilde{k}_\theta - h_\theta\| = \sqrt{\int_0^1 \tilde{k}_\theta^2(x)dx - \sum_{n=1}^{N(\theta)} k_n^2(\theta)}$, $\sigma = 1$, $\beta = 1/10$ and requiring that

$\gamma\|\tilde{k}_\theta - h_\theta\|\tilde{L}_\theta \leq 1 - \frac{1}{12M}$. Fig. 2 shows the evolution of $\|u[t]\|$ for the closed-loop system (2.1), (2.2), (2.3) with (2.9), (2.10), (2.13), (2.19), (2.20), (4.2), (4.3), (4.9), (4.10), (4.11) with the following choices for the tunable parameters and the initial condition:

$$T = 0.05, \ a = 1, \ M = 10, \ \tilde{N} = 1, \ u_0(x) = \sqrt{2}\sin(\pi x) + x^2 - x^3, \ x \in [0,1], \ \hat{\theta}_0 = 0.1 \quad (4.12)$$

Fig. 2 shows the difference between the response of the closed-loop system with the regulation-triggered adaptive scheme and the closed-loop system with the nominal feedback (for known $\theta$). The event trigger becomes active at $t = T$ and the exact value of the parameter $\theta = 11$ is found at that time. Later than this time, the adaptive controller coincides with the nominal feedback.

Finally, for comparison purposes we also simulated the adaptive controller given in [32,34] with the passive identifier, namely, the controller

$$\frac{\partial \hat{u}}{\partial t}(t,x) = \frac{\partial^2 \hat{u}}{\partial x^2}(t,x) + \hat{\theta}(t)u(t,x) + \gamma^2 \|u[t]\|^2 (u(t,x) - \hat{u}(t,x)), \ t > 0, \ x \in (0,1) \quad (4.13)$$

$$\hat{u}(t,0) = \hat{u}(t,1) - U(t) = 0, \text{ for } t > 0 \quad (4.14)$$

$$U(t) = -\hat{\theta}(t)\int_0^1 s \frac{I_1\left(\sqrt{\hat{\theta}(t)(1-s^2)}\right)}{\sqrt{\hat{\theta}(t)(1-s^2)}} u(t,s)ds \quad (4.15)$$

$$\frac{d\hat{\theta}}{dt}(t) = \gamma \int_0^1 (u(t,x) - \hat{u}(t,x))u(t,x)dx \quad (4.16)$$

with the following parameter values and same initial conditions as above:

$$\gamma = 100, \ u_0(x) = \hat{u}(0,x) = \sqrt{2}\sin(\pi x) + x^2 - x^3, \ x \in [0,1], \ \hat{\theta}_0 = \hat{\theta}(0) = 0.1 \quad (4.17)$$

The results are shown in Fig.3 and Fig. 4. The overshoot of $L^2$ norm for the adaptive controller with the passive identifier is less than the overshoot of the $L^2$ norm for the regulation-triggered adaptive controller but the convergence of the state to zero is much slower. Moreover, the estimation of the value of the parameter is not accurate for the passive identifier: the passive identifier heavily overestimates the value of the unknown parameter.

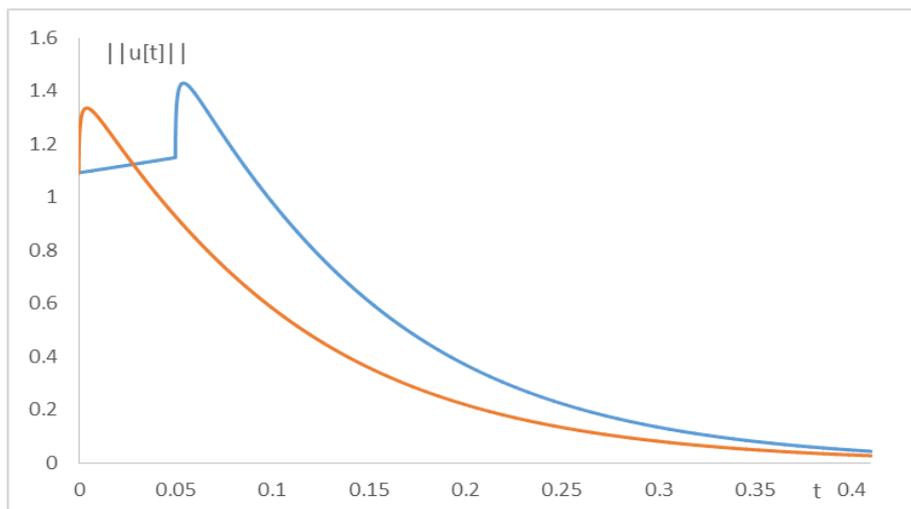

**Fig. 2:** The blue curve shows the evolution of $\|u[t]\|$ for the closed-loop system (2.1), (2.2), (2.3) with (2.9), (2.10), (2.13), (2.19), (2.20), (4.2), (4.3), (4.9), (4.10), (4.11) and (4.12). The red curve shows the evolution of $\|u[t]\|$ for the closed-loop system (2.1), (2.2), (2.3) with (2.4), (2.5), (4.8) and same initial condition.

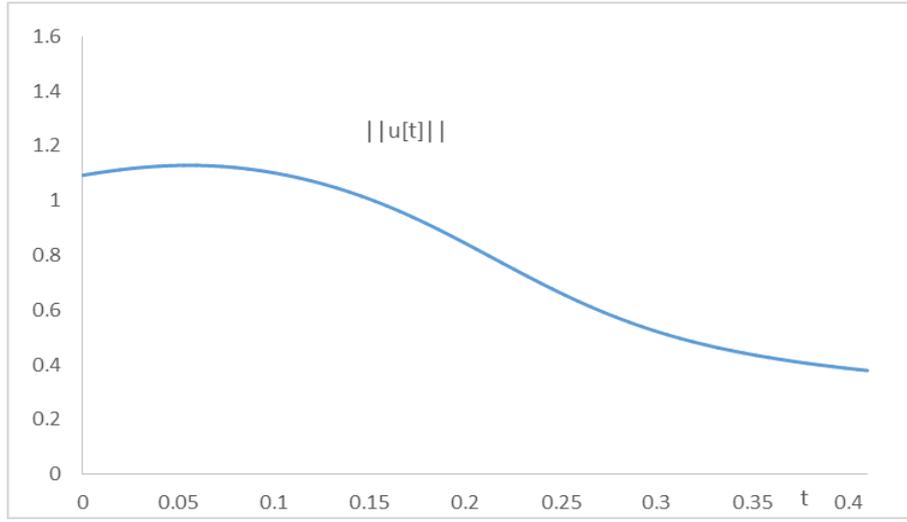

**Fig. 3:** The evolution of $\|u[t]\|$ for the closed-loop system (2.1), (2.2), (2.3) with (4.13), (4.14), (4.15), (4.16) and (4.17).

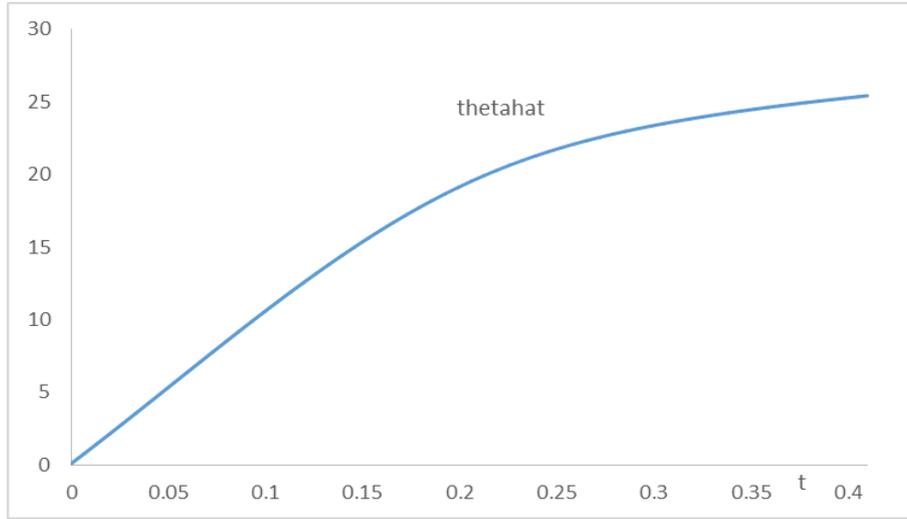

**Fig. 4:** The evolution of the parameter estimation $\hat{\theta}(t)$ for the closed-loop system (2.1), (2.2), (2.3) with (4.13), (4.14), (4.15), (4.16) and (4.17).

## 5. Proof of Main Result

For the proof of Theorem 2.2, we need first to develop certain auxiliary results. The first auxiliary results deals with the parabolic equation

$$\frac{\partial u}{\partial t}(t,x) = p\frac{\partial^2 u}{\partial x^2}(t,x), \text{ for } t>0, \ x\in(0,1) \tag{5.1}$$

where $p>0$ is a constant, with boundary conditions

$$u(t,0) = 0, \text{ for } t>0 \tag{5.2}$$
$$u(t,1) = U(t), \text{ for } t>0 \tag{5.3}$$

where $U(t)$ is the control input, which is given by a feedback law of the form

$$U(t) = \int_0^1 k(s)u(t,s)ds \tag{5.4}$$

where

$$k(x) = \sum_{n=1}^{N} k_n \phi_n(x), \quad x \in [0,1] \tag{5.5}$$

the functions $\phi_n$ are defined by (2.6) and $k_n$, $n = 1, 2, ..., N$ are constants.

**Theorem 5.1:** *Consider system (5.1), (5.2), (5.3), where $p > 0$ is a constant. Then for every integer $N \geq 1$ and for any $(k_1, k_2, ..., k_N) \in \Re^N$ there exist constants $Q, \sigma > 0$ such that for every $u_0 \in L^2(0,1)$, the initial boundary value problem (5.1), (5.2), (5.3), (5.4) with initial condition $u[0] = u_0$, where $k$ is defined by (5.5), (2.6), has a unique solution $u \in C^0(\Re_+; L^2(0,1)) \cap C^1((0, +\infty) \times [0,1])$ with $u[t] \in C^2([0,1])$ for $t > 0$ which also satisfies the following estimate*

$$\|u[t]\| \leq Q \exp(\sigma t) \|u_0\|, \quad t \geq 0 \tag{5.6}$$

The proof of the Theorem 5.1 is based on some technical lemmas, whose proofs can be found at the Appendix.

**Lemma 5.2:** *Let $k \in C^2([0,1])$ be a given function. There exists a constant $\sigma > 0$ such that for every solution $w \in C^0(\Re_+; L^2(0,1)) \cap C^1((0, +\infty) \times [0,1])$ with $w[t] \in C^2([0,1])$ for $t > 0$ of the following initial-boundary value problem*

$$\frac{\partial w}{\partial t}(t,x) = p \frac{\partial^2 w}{\partial x^2}(t,x) + a\left(k''(x) + (x + \gamma k(x))\|k'\|^2\right) \int_0^1 k(s) w(t,s) ds - p(x + \gamma k(x)) \int_0^1 k''(s) w(t,s) ds,$$

$$\text{for } t > 0, \; x \in (0,1) \tag{5.7}$$

*where $p > 0$, $a, \gamma \in \Re$ are constants with*

$$w[0] = w_0 \in L^2(0,1) \tag{5.8}$$

$$w(t,0) = w(t,1) = 0, \text{ for } t > 0 \tag{5.9}$$

*the following estimate holds*

$$\|w[t]\| \leq \exp(\sigma t) \|w_0\|, \text{ for all } t \geq 0 \tag{5.10}$$

**Lemma 5.3:** *Let $N \geq 1$ be an integer and let $(k_1, k_2, ..., k_N) \in \Re^N$ be given. There exist constants $Q, \sigma > 0$ such that, for every solution $u \in C^0(\Re_+; L^2(0,1)) \cap C^1((0, +\infty) \times [0,1])$ with $u[t] \in C^2([0,1])$ for $t > 0$ of the initial-boundary value problem (5.1), (5.2), (5.3), (5.4), with initial condition $u[0] = u_0 \in L^2(0,1)$, where $k \in C^2([0,1])$ is defined by (5.5), (2.6), estimate (5.6) holds.*

**Lemma 5.4:** *Let $N \geq 1$ be an integer and let $(k_1, k_2, ..., k_N) \in \Re^N$ be given. For every $u_0 \in L^2(0,1)$, there exists at most one solution $u \in C^0(\Re_+; L^2(0,1)) \cap C^1((0, +\infty) \times [0,1])$ with $u[t] \in C^2([0,1])$ for $t > 0$ of the initial boundary value problem (5.1), (5.2), (5.3), (5.4), with initial condition $u[0] = u_0 \in L^2(0,1)$, where $k \in C^2([0,1])$ is defined by (5.5), (2.6).*

**Lemma 5.5:** *Let $p > 0$, $a, \gamma \in \Re$ be given constants, let $N \geq 1$ be an integer and let $(k_1, k_2, ..., k_N) \in \Re^N$ be given. Define the vector*

$$r := \pi^2 [k_1(ab + p), ..., k_N(ab + pN^2)] \tag{5.11}$$

*and the matrix $A = \{A_{i,j} : i, j = 1, ... N\} \in \Re^{N \times N}$ by means of the formula*

$$A_{i,j} := k_i k_j \pi^2 \left(a\gamma b + p\gamma j^2 - i^2 a\right) + L_i k_j \pi^2 \left(ab + pj^2\right) - i^2 \pi^2 p \delta_{i,j} \tag{5.12}$$

*where $\delta_{i,j} = 1$ if $i = j$ and $\delta_{i,j} = 0$ if $i \neq j$,*

$$b := \sum_{n=1}^{N} n^2 k_n^2 \tag{5.13}$$

$$L_n = -\frac{(-1)^n \sqrt{2}}{n\pi}, \quad n = 1, 2, \dots \tag{5.14}$$

*For any $w_0 \in L^2(0,1)$, $T > 0$, the function $w : [0,T] \times [0,1] \to \Re$ defined by the following formula for $x \in [0,1]$, $t \in (0,T]$*

$$w(t,x) = \sum_{n=1}^{N} \phi_n(x) e_n \exp(At) \xi + \sum_{n=N+1}^{\infty} \exp\left(-pn^2\pi^2 t\right) c_n \phi_n(x)$$
$$+ \sum_{n=N+1}^{\infty} L_n \phi_n(x) \left( \int_0^t \exp\left(-pn^2\pi^2(t-s)\right) r \exp(As) \xi \, ds \right) \tag{5.15}$$

*and*

$$w(0,x) = w_0(x), \text{ for } x \in [0,1] \tag{5.16}$$

*where $\phi_n$ are given by (2.6), $e_n \in \Re^N$, $n = 1, \dots N$, and $e_1 = (1, 0 \dots, 0), \dots, e_N = (0, \dots 0, 1)$,*

$$c_n := \int_0^1 w_0(x) \phi_n(x) dx, \quad n = 1, 2, \dots \tag{5.17}$$

$$\xi := \begin{bmatrix} c_1 \\ \vdots \\ c_N \end{bmatrix} \tag{5.18}$$

*is of class $C^0((0,T] \times [0,1])$. Moreover, the mapping $[0,T] \ni t \to w[t] \in L^2(0,1)$ is continuous and $w[t] \in C^1([0,1])$ for every $t \in (0,T]$.*

Using Lemma 5.2, Lemma 5.3, Lemma 5.4 and Lemma 5.5, we are in a position to prove Theorem 5.1.

**Proof of Theorem 5.1:** Uniqueness follows from Lemma 5.4. Let $\gamma \in \Re$ be a constant with $\int_0^1 xk(x)dx + \gamma \|k\|^2 \neq 1$ (notice that for every $k \in C^2([0,1])$ there are infinite values of $\gamma \in \Re$ for which this relation holds) and define

$$\beta := \left( 1 - \int_0^1 sk(s)ds - \gamma \|k\|^2 \right)^{-1} \tag{5.19}$$

$$a := \beta p \gamma \tag{5.20}$$

Consider the following Fredholm transformation

$$w(t,x) = u(t,x) - (x + \gamma k(x)) \int_0^1 k(s) u(t,s) ds, \text{ for } t \geq 0, \ x \in [0,1] \tag{5.21}$$

and its inverse

$$u(t,x) = w(t,x) + \beta (x + \gamma k(x)) \int_0^1 k(s) w(t,s) ds, \text{ for } t \geq 0, \ x \in [0,1] \tag{5.22}$$

Using (5.21), (5.22), (5.1), (5.2), (5.3), (5.4), (5.5), we guarantee that the initial-boundary value problem (5.1), (5.2), (5.3), (5.4) with initial condition $u[0] = u_0$ is transformed to the initial-boundary value problem (5.7), (5.8), (5.9) with

$$w_0(x) = u_0(x) - (x + \gamma k(x)) \int_0^1 k(s) u_0(s) ds, \text{ for } x \in [0,1] \tag{5.23}$$

Therefore, it suffices to show that for every $T > 0$ there exists a function $w \in C^0([0,T]; L^2(0,1)) \cap C^1((0,T) \times [0,1])$ with $w[t] \in C^2([0,1])$ for $t \in (0,T]$ which satisfies the following equations:

$$\frac{\partial w}{\partial t}(t,x) = p \frac{\partial^2 w}{\partial x^2}(t,x) + a\left(k''(x) + (x + \gamma k(x)) \|k'\|^2\right) \int_0^1 k(s) w(t,s) ds - p(x + \gamma k(x)) \int_0^1 k''(s) w(t,s) ds,$$

$$\text{for } t \in (0,T], \ x \in (0,1) \tag{5.24}$$

$$w[0] = w_0 \in L^2(0,1) \tag{5.25}$$

$$w(t,0) = w(t,1) = 0, \text{ for } t \in (0,T] \tag{5.26}$$

with $w_0$ defined by (5.23). The fact that (5.6) holds is a direct consequence of Lemma 5.3. In what follows $T > 0$ is arbitrary.

We next claim that the function $w:[0,T] \times [0,1] \to \Re$ defined by (5.15), (5.16) is the required function. First, we notice that by virtue of Lemma 5.5, it holds that $w \in C^0([0,T]; L^2(0,1)) \cap C^0((0,T] \times [0,1])$ with $w[t] \in C^1([0,1])$ for $t \in (0,T]$. Next, we show that the function $w:[0,T] \times [0,1] \to \Re$ defined by (5.15), (5.16) satisfies $w[t] \in C^2([0,1])$ for every $t \in (0,T]$ and that $\frac{\partial w}{\partial t}(t,x)$ is a continuous function on $(0,T] \times [0,1]$. According to definitions (2.6) and (5.14), the following equality holds

$$\frac{x - x^3}{6} = \sum_{n=1}^{\infty} \frac{L_n \phi_n(x)}{n^2 \pi^2}, \text{ for } x \in [0,1] \tag{5.27}$$

Using integration by parts in (5.15), and taking into account (5.27), we obtain the following formula

$$w(t,x) = \sum_{n=1}^{N} \phi_n(x) e_n \exp(At) \xi + \sum_{n=N+1}^{\infty} \exp(-pn^2\pi^2 t) c_n \phi_n(x)$$

$$+ r \exp(At) \xi \frac{x - x^3}{6p} - r\xi \sum_{n=N+1}^{\infty} L_n \phi_n(x) \frac{\exp(-pn^2\pi^2 t)}{pn^2\pi^2} - r\exp(At)\xi \sum_{n=1}^{N} \frac{L_n \phi_n(x)}{pn^2\pi^2}$$

$$- \sum_{n=N+1}^{\infty} \frac{L_n \phi_n(x)}{pn^2\pi^2} \left( \int_0^t \exp(-pn^2\pi^2(t-s)) rA \exp(As) \xi\, ds \right)$$

$$\text{for } x \in [0,1], \ t \in (0,T] \tag{5.28}$$

Differentiating formally (5.28) twice term by term with respect to $x \in [0,1]$ and once with respect to $t \in (0,T]$ and taking into account definition (2.6) and (5.27), we get

$$\frac{\partial^2 w}{\partial x^2}(t,x) = -\sum_{n=1}^{N} n^2 \pi^2 \phi_n(x) e_n \exp(At) \xi - \sum_{n=N+1}^{\infty} n^2 \pi^2 \exp(-pn^2\pi^2 t) c_n \phi_n(x)$$

$$- r \exp(At) \xi \frac{x}{p} + \frac{r\xi}{p} \sum_{n=N+1}^{\infty} L_n \phi_n(x) \exp(-pn^2\pi^2 t) + \frac{r \exp(At) \xi}{p} \sum_{n=1}^{N} L_n \phi_n(x) \tag{5.29}$$

$$+ \frac{1}{p} \sum_{n=N+1}^{\infty} L_n \phi_n(x) \left( \int_0^t \exp(-pn^2\pi^2(t-s)) rA \exp(As) \xi\, ds \right)$$

$$\frac{\partial w}{\partial t}(t,x) = \sum_{n=1}^{N} \phi_n(x) e_n A \exp(At) \xi - \sum_{n=N+1}^{\infty} pn^2\pi^2 \exp(-pn^2\pi^2 t) c_n \phi_n(x)$$

$$+ r\xi \sum_{n=N+1}^{\infty} L_n \phi_n(x) \exp(-pn^2\pi^2 t) + \sum_{n=N+1}^{\infty} L_n \phi_n(x) \left( \int_0^t \exp(-pn^2\pi^2(t-s)) rA \exp(As) \xi\, ds \right) \tag{5.30}$$

At this point, we show that for any $w_0 \in L^2(0,1)$, $T > 0$ the function $w:[0,T] \times [0,1] \to \Re$ defined by (5.15), (5.16) satisfies (5.29) and simultaneously, $\frac{\partial w}{\partial t}(t,x)$ is a continuous function on $(0,T] \times [0,1]$ and satisfies (5.30). Combining (5.14), (5.17), the fact that $w_0 \in L^2(0,1)$ and the fact that $\|\phi_n\| = 1$ for $n = 1, 2, \ldots$ $x^a \exp(-x) \leq a^a \exp(-a)$, for any $x \geq 0$, $a > 0$ and applying the Cauchy-Schwartz inequality, we get for all $t_0 \in (0,T]$ and $t \in [t_0, T]$:

$$\left|\left(n^2\pi^2|c_n| + |L_n|\right)\exp\left(-pn^2\pi^2 t\right)\right| \leq \frac{4\exp(-2)}{p^2 n^2 \pi^2 t_0^2}\|w_0\| + \frac{\sqrt{2}\exp(-1)}{pn^3\pi^3 t_0}, \text{ for } n = 1, 2, \ldots \quad (5.31)$$

Moreover, combining definition (5.14) along with the fact that $\|\phi_n\| = 1$ for $n = 1, 2, \ldots$, $|\exp(At)| \leq \exp(|A|t)$ for all $t \geq 0$, we obtain for all $t_0 \in (0,T]$ and $t \in [t_0, T]$:

$$\left|L_n \int_0^t \exp\left(-pn^2\pi^2(t-s)\right) rA\exp(As)\xi\, ds\right| \leq |r|\|A\|\|\xi\|\exp(|A|T)\frac{\sqrt{2}}{pn^3\pi^3}, \text{ for } n = 1, 2, \ldots \quad (5.32)$$

Inequalities (5.31), (5.32) indicate that the infinite series appearing on the right hand side of (5.29), (5.30) are uniformly and absolutely convergent on $[t_0, T] \times [0,1]$, for every $t_0 \in (0,T]$. Therefore $w[t] \in C^2([0,1])$ for every $t \in (0,T]$ and $\frac{\partial w}{\partial t}(t,x)$ is a continuous function on $(0,T] \times [0,1]$.

Using (5.29), (5.30), we obtain for all $t \in (0,T]$, $x \in (0,1)$

$$\frac{\partial w}{\partial t}(t,x) = p\frac{\partial^2 w}{\partial x^2}(t,x) + \sum_{n=1}^{N} \phi_n(x) e_n A\exp(At)\xi$$
$$+ r\exp(At)\xi x - r\exp(At)\xi \sum_{n=1}^{N} L_n \phi_n(x) + p\sum_{n=1}^{N} n^2\pi^2 \phi_n(x) e_n \exp(At)\xi \quad (5.33)$$

Notice that definitions (5.12), (5.11) imply that $e_n A = -n^2\pi^2 p e_n + L_n r + \gamma k_n r - an^2\pi^2 k_n \bar{k}$ for all $n = 1, \ldots, N$, where $\bar{k} := [k_1, \ldots, k_N]$. Combining (5.33) with definitions (5.5), (2.6) and the previous equalities, we get for all $t \in (0,T]$, $x \in (0,1)$

$$\frac{\partial w}{\partial t}(t,x) = p\frac{\partial^2 w}{\partial x^2}(t,x) + r\exp(At)\xi\left(x + \gamma k(x)\right) + a\bar{k}\exp(At)\xi k''(x) \quad (5.34)$$

Using (5.5), (2.6), (5.11) and (5.15) we get for all $t \in (0,T]$, $x \in (0,1)$

$$b\pi^2 = \|k'\|^2$$
$$\bar{k}\exp(At)\xi = \int_0^1 k(s)w(t,s)ds \quad (5.35)$$
$$r\exp(At)\xi = ab\pi^2\int_0^1 k(s)w(t,s)ds - p\int_0^1 k''(s)w(t,s)ds$$

It follows from (5.34), (5.35) that (5.24) holds. Equalities (5.25), (5.26) also hold.

Therefore, we conclude that the function $w:[0,T] \times [0,1] \to \Re$ defined by (5.15), (5.16), is a function of class $C^0([0,T]; L^2(0,1)) \cap C^1((0,T] \times [0,1])$ with $w[t] \in C^2([0,1])$ for $t \in (0,T]$ which satisfies (5.24), (5.25), (5.26). Consequently, according to the inverse transformation (5.22) and definition (5.5),(2.6) which implies that $k \in C^2([0,1])$, we conclude that $u \in C^0(\Re_+; L^2(0,1)) \cap C^1((0,+\infty) \times [0,1])$ with $u[t] \in C^2([0,1])$ for $t > 0$ and satisfies the initial-boundary value problem (5.1), (5.2), (5.3), (5.4) where $k \in C^2([0,1])$ is defined by (5.5), (2.6). The proof is complete. ◁

Our second auxiliary result is the following corollary.

**Corollary 5.6:** *Let $N \geq 1$ be an integer and let $(k_1, k_2,..., k_N) \in \Re^N$ be a given vector. Then there exist constants $Q, \sigma > 0$ such that for every $\tilde{u}_0 \in L^2(0,1)$, the initial-boundary value problem*

$$\frac{\partial \tilde{u}}{\partial t}(t,x) = p\frac{\partial^2 \tilde{u}}{\partial x^2}(t,x) + \theta \tilde{u}(t,x), \quad t > 0, \quad x \in (0,1) \tag{5.36}$$

*where $p > 0$, $\theta \in \Re$ are constants, with boundary conditions*

$$\tilde{u}(t,0) = 0, \text{ for } t > 0 \tag{5.37}$$

$$\tilde{u}(t,1) = \tilde{U}(t), \text{ for } t > 0 \tag{5.38}$$

*and $\tilde{U}(t)$ being given by the formula*

$$\tilde{U}(t) = \int_0^1 k(s)\tilde{u}(t,s)ds, \text{ for } t > 0 \tag{5.39}$$

*where $k:[0,1] \to \Re$ is defined by (5.5), (2.6), with initial condition $\tilde{u}[0] = \tilde{u}_0 \in L^2(0,1)$ has a unique solution $\tilde{u} \in C^0(\Re_+; L^2(0,1)) \cap C^1((0,+\infty) \times [0,1])$ with $\tilde{u}[t] \in C^2([0,1])$ for $t > 0$, which also satisfies estimate (5.6).*

**Proof:** Define $\tilde{u}(t,x) = \exp(\theta t)u(t,x)$, where $u \in C^0(\Re_+; L^2(0,1)) \cap C^1((0,+\infty) \times [0,1])$ with $u[t] \in C^2([0,1])$ for $t > 0$ is the unique solution of the initial-boundary value problem (5.1), (5.2), (5.3), (5.4) with initial condition $u[0] = \tilde{u}_0$. The rest is a direct consequence of Theorem 5.1. ◁

Our third auxiliary result deals with the case that the input is identically zero.

**Lemma 5.7:** *Let $N \geq 1$ be an integer let $p > 0$, $\theta \in \Re$, $T > 0$ be constants and let $(k_1, k_2,..., k_N) \in \Re^N$ be a given vector. Suppose that the solution of the initial-boundary value problem (5.36), (5.37), (5.38), (5.39), where $k:[0,1] \to \Re$ is defined by (5.5), (2.6), with initial condition $\tilde{u}[0] = \tilde{u}_0 \in L^2(0,1)$, satisfies $\tilde{U}(t) = 0$ for $t \in [0,T)$. Then $k_n \int_0^1 \sin(n\pi x)\tilde{u}(t,x)dx = 0$ for all $t \in [0,T)$, $n = 1,...,N$.*

**Proof:** Define:

$$a_n(t) := \int_0^1 \sin(n\pi x)\tilde{u}(t,x)dx, \text{ for } n = 1,...,N \tag{5.40}$$

Using (5.40), (5.36), (5.37), (5.38), the fact that $\tilde{U}(t) = 0$ for $t \in [0,T)$ and integration by parts we get $\dot{a}_n(t) = (\theta - n^2\pi^2 p)a_n(t)$, for $n = 1,...,N$ and $t \in (0,T)$. Using the fact that $\tilde{u} \in C^0(\Re_+; L^2(0,1))$ (which implies that the mappings $t \to a_n(t)$ are continuous) and integrating the differential equations, we obtain:

$$a_n(t) = \exp\left((\theta - n^2\pi^2 p)t\right)a_n(0), \text{ for } n = 1,...,N \tag{5.41}$$

Using (5.5), (2.6), (5.40), the fact that $\tilde{U}(t) = 0$ for $t \in [0,T)$ and (5.41), we get:

$$\sum_{n=1}^N k_n (n^2\pi^2)^m a_n(0) = 0, \text{ for } m = 0,...,N-1 \tag{5.42}$$

Equations (5.42) can be written as $Ax=0$, where $A$ is the transpose of the invertible Vandermonde matrix $\begin{bmatrix} 1 & (1^2\pi^2) & \cdots & (1^2\pi^2)^{N-1} \\ 1 & (2^2\pi^2) & \cdots & (2^2\pi^2)^{N-1} \\ \vdots & \vdots & & \vdots \\ 1 & (N^2\pi^2) & \cdots & (N^2\pi^2)^{N-1} \end{bmatrix}$ and $x = \begin{bmatrix} k_1 a_1(0) \\ k_2 a_2(0) \\ \vdots \\ k_N a_N(0) \end{bmatrix}$. Since $A$ is invertible we get $x=0$ and $k_n a_n(0) = 0$ for all $n \in I$. The conclusion is a direct consequence of (5.41) and definition (5.40). ◁

We are now ready to provide the proof of Theorem 2.2.

**Proof of Theorem 2.2:** By a solution $(u[t], \hat{\theta}(t), \hat{c}(t))$ of (2.1), (2.2), (2.3), (2.8), (2.9), (2.10), (2.12), (2.13), (2.19), (2.28) on an interval $[0, \bar{t}]$ with $\bar{t} > 0$, we mean mappings $\hat{\theta}, \hat{c} : [0, \bar{t}] \to \Re$, $u \in C^0([0,\bar{t}]; L^2(0,1))$ satisfying $u \in C^1(I_{\bar{t}} \times [0,1])$, $u[t] \in C^2([0,1])$ for $t \in (0, \bar{t}]$, where $I_{\bar{t}} = (\Re_+ \setminus \{\tau_i \geq 0, i = 0,1,2,...\}) \cap [0,\bar{t}]$, which also satisfies (2.8), (2.9), (2.10), (2.12), (2.13), (2.19), (2.28) for all $t \in [0, \bar{t}]$ and for all $i \in Z_+$ with $\tau_i \leq \bar{t}$ and

$$\frac{\partial u}{\partial t}(t,x) = p \frac{\partial^2 u}{\partial x^2}(t,x) + \theta u(t,x), \text{ for all } (t,x) \in I_{\bar{t}} \times (0,1) \tag{5.43}$$

$$u(t,0) = u(t,1) - cU(t) = 0, \text{ for all } t \in I_{\bar{t}} \tag{5.44}$$

Similarly, we may define the solution $(u[t], \hat{\theta}(t), \hat{c}(t))$ of (2.1), (2.2), (2.3), (2.8), (2.9), (2.10), (2.12), (2.13), (2.19), (2.28) on an interval $[0, \bar{t})$ with $\bar{t} > 0$ (replace all inequalities that involve $\bar{t}$ above with strict inequalities and all replace the interval $[0, \bar{t}]$ by $[0, \bar{t})$ in all relations).

**Claim 1:** *If a solution $(u[t], \hat{\theta}(t), \hat{c}(t))$ of (2.1), (2.2), (2.3), (2.8), (2.9), (2.10), (2.12), (2.13), (2.19), (2.28) is defined on $t \in [0, \tau_i]$ for certain $i \in Z_+$, then the solution is defined on $t \in [0, \tau_{i+1}]$. Moreover, it holds that*

$$\|u[t]\| \leq R_{\hat{\theta}(\tau_i)}(1+a)\|u[\tau_i]\|, \text{ for all } t \in [\tau_i, \tau_{i+1}] \tag{5.45}$$

**Proof of Claim 1:** Consider the initial boundary value problem

$$\frac{\partial w}{\partial t}(t,x) = p \frac{\partial^2 w}{\partial x^2}(t,x) + \theta w(t,x), \ t > 0, \ x \in (0,1) \tag{5.46}$$

$$w(t,0) = w(t,1) - \frac{c}{\hat{c}(\tau_i)} \int_0^1 k(\hat{\theta}(\tau_i), x) w(t,x) dx = 0, \text{ for } t > 0 \tag{5.47}$$

$$w[0] = u[\tau_i], \text{ for } t > 0 \tag{5.48}$$

where $k(\theta, x)$ is given by (2.5). By virtue of Corollary 5.6, there exists a unique solution $w \in C^0(\Re_+; L^2(0,1)) \cap C^1((0,+\infty) \times [0,1])$ with $w[t] \in C^2([0,1])$ for $t > 0$. Define $\tau_{i+1} > \tau_i$ by means of (2.10), where

$$r_i := \tau_i + \inf\left\{s > 0 : \|w[s]\| = R_{\hat{\theta}(\tau_i)}(1+a)\|w[0]\|\right\}, \text{ for } w[0] \neq 0 \tag{5.49}$$

$$r_i := \tau_i + T, \text{ for } w[0] = 0 \tag{5.50}$$

Next, define

$$u[\tau_i + s] = w[s], \text{ for } s \in (0, \tau_{i+1} - \tau_i] \tag{5.51}$$

and $(\hat{\theta}(t), \hat{c}(t))$ for $t \in (\tau_i, \tau_{i+1}]$ by means of (2.9), (2.19) and (2.28).

It can be verified that the mapping $(u[t], \hat{\theta}(t), \hat{c}(t))$ is a solution of (2.1), (2.2), (2.3), (2.8), (2.9), (2.10), (2.12), (2.13), (2.19), (2.28) defined on $t \in [0, \tau_{i+1}]$. Moreover, using (5.51), (5.48) and (5.49), (5.50), it follows that (5.45) holds. The claim is proved. ◁

It follows from Claim 1 (and a step-by-step construction) that for every $u_0 \in L^2(0,1)$, $\hat{\theta}_0 \in \Re$, $\hat{c}_0 \in \Re$ the solution of (2.1), (2.2), (2.3), (2.8), (2.9), (2.10), (2.12), (2.13), (2.19), (2.28) with initial conditions $u[0] = u_0$, $\hat{\theta}(0) = \hat{\theta}_0$, $\hat{c}(0) = \hat{c}_0$ is unique and is defined on $\left[0, \lim_{i \to +\infty}(\tau_i)\right)$.

The following two claims are direct consequences of Lemma 2.1.

**Claim 2:** *If a solution $(u[t], \hat{\theta}(t), \hat{c}(t))$ of (2.1), (2.2), (2.3), (2.8), (2.9), (2.10), (2.12), (2.13), (2.19), (2.28) satisfies $\hat{\theta}(\tau_i) = \theta$ (or $\hat{c}(\tau_i) = c$) for certain $i \in Z_+$, then the solution satisfies $\hat{\theta}(t) = \theta$ (or $\hat{c}(t) = c$) for all $t \in \left[\tau_i, \lim_{k \to +\infty}(\tau_k)\right)$.*

**Claim 3:** *If $u_0 \neq 0$ then $\hat{\theta}(\tau_1) = \theta$.*

We next show the following claim.

**Claim 4:** *If $u_0 \neq 0$ and $\hat{c}(\tau_2) \neq c$ then $U(t) = 0$ for $t \in \left[0, \lim_{k \to +\infty}(\tau_k)\right)$ and (2.4) holds for $t \in \left[\tau_1, \lim_{k \to +\infty}(\tau_k)\right)$.*

**Proof of Claim 4:** It follows from Claim 2 that $\hat{c}(\tau_l) \neq c$ for $l = 0, 1$. By virtue of Lemma 2.1, we have that $U(t) = 0$ for $t \in [0, \tau_2)$. Lemma 5.7, Claim 2 and Claim 3 in conjunction with (2.5), (2.6), (2.8), (2.9) imply that

$$k_n(\theta) \int_0^1 \sin(n\pi x) u(t,x) dx = 0 \text{ for all } t \in [\tau_1, \tau_2) \text{ and } n = 1, \ldots, N(\theta) \qquad (5.52)$$

Notice that the unique solution of (2.1), (2.2), (2.3), (2.8), (2.9), (2.10), (2.12), (2.13), (2.19), (2.28) on the interval $t \in \left[\tau_1, \lim_{k \to +\infty}(\tau_k)\right)$ is given by the formula $u(t,x) = 2 \sum_{n=1}^{\infty} \exp\left((\theta - n^2 \pi^2 p)(t - \tau_1)\right) \sin(n\pi x) \int_0^1 \sin(n\pi z) u(\tau_1, z) dz$. Therefore, we conclude from (5.52), (2.5), (2.6), (2.8), (2.9) that $U(t) = 0$ for $t \in \left[\tau_1, \lim_{k \to +\infty}(\tau_k)\right)$. Again, by virtue of (5.52), it follows that (2.4) holds for $t \in \left[\tau_1, \lim_{k \to +\infty}(\tau_k)\right)$. The claim is proved. ◁

Therefore, Claim 2 and Claim 4 implies that (2.4) holds for all $t \in \left[\tau_2, \lim_{k \to +\infty}(\tau_k)\right)$, when $u_0 \neq 0$. The following claim specifies the times of the events.

**Claim 5:** *If $u_0 \neq 0$, then the solution is defined for all $t \geq 0$ and satisfies $\tau_j = \tau_1 + (j-2)T$ for $j \geq 2$.*

**Proof of Claim 5:** We prove by induction that $\tau_{i+1} = \tau_i + T$ for $i \geq 2$. Let $i \geq 2$ be an integer. Notice that Claim 2 and Claim 4 imply that $\hat{\theta}(t) = \theta$ and that (2.4) holds for all $t \in [\tau_i, \tau_{i+1})$. First assume that $u[\tau_i] \neq 0$. By virtue of (2.7) and since (2.4) holds, we get $\|u[t]\| \leq R_{\hat{\theta}(\tau_i)} \|u[\tau_i]\|$, for all $t \in [\tau_i, \tau_{i+1})$. It follows that $\|u[t]\| \leq R_{\hat{\theta}(\tau_i)} \|u[\tau_i]\| < R_{\hat{\theta}(\tau_i)}(1+a)\|u[\tau_i]\|$ for all $t \in [\tau_i, \tau_{i+1})$. Therefore, we get from (2.10), (2.12) that $\tau_{i+1} = \tau_i + T$. The same conclusion follows from (2.10), (2.13) if $u[\tau_i] = 0$. The claim is proved. ◁

By virtue of Claim 5, we conclude that $\lim_{i \to +\infty}(\tau_i) = +\infty$ when $u_0 \neq 0$. Moreover, it follows from Claim 5 and (2.7) that there exist constants $R_\theta \geq 1$, $\omega_\theta > 0$, such that the following estimate holds when $u_0 \neq 0$:

$$\|u[t]\| \leq R_\theta \exp(-\omega_\theta (t-\tau_2))\|u[\tau_2]\|, \text{ for } t \geq \tau_2 \tag{5.53}$$

Using Corollary 5.6 and Lemma 2.1 (which implies that $\hat{c}(\tau_1) = c$ or $\hat{c}(\tau_1) = \hat{c}_0$), it follows that there exist constants $Q(\hat{\theta}_0, \hat{c}_0, \theta, c), \sigma(\hat{\theta}_0, \hat{c}, \theta, c) > 0$ depending on $\hat{\theta}_0 \in \Re$ such that

$$\|u[t]\| \leq Q(\hat{\theta}_0, \hat{c}_0, \theta, c)\exp(\sigma(\hat{\theta}_0, \hat{c}_0, \theta, c)t)\|u_0\|, \text{ for } t \in [0, \tau_2] \tag{5.54}$$

Inequality (5.53) in conjunction with (5.54) for $t = \tau_1$ implies the following estimate when $u_0 \neq 0$:

$$\|u[t]\| \leq R_\theta \exp((\omega_\theta + \sigma(\hat{\theta}_0, \hat{c}_0, \theta, c))\tau_2)Q(\hat{\theta}_0, \hat{c}_0, \theta, c)\exp(-\omega_\theta t)\|u_0\|, \text{ for } t \geq 0 \tag{5.55}$$

Estimate (5.55) holds even if $u_0 = 0$, because in this case the solution satisfies $\|u[t]\| = 0$ for all $t \geq 0$. Setting $M_{\theta,\hat{\theta}_0,c,\hat{c}_0} = R_\theta \exp((\omega_\theta + \sigma(\hat{\theta}_0, \hat{c}_0, \theta, c))\tau_2)Q(\hat{\theta}_0, \hat{c}_0, \theta, c) > 0$, we obtain directly from (5.55) that the estimate $\|u[t]\| \leq M_{\theta,\hat{\theta}_0} \exp(-\omega_\theta t)\|u_0\|$ holds for all $t \geq 0$.

Finally, if $u_0 \neq 0$, using Claim 2 and Claim 3, we obtain the equality $\hat{\theta}(t) = \theta$ for all $t \geq T$. Moreover, if there exists $t \geq 0$ with $U(t) \neq 0$ then Claim 2 and Claim 4 imply that $\hat{c}(t) = c$ for all $t \geq \tau_2$ (the case $u_0 = 0$ is excluded because in this case $U(t) \equiv 0$). The proof is complete. ◁

## 6. Concluding Remarks

We have developed a novel adaptive boundary control scheme for parabolic PDEs, which guarantees exponential convergence of the state to zero in the $L^2$ norm and finite-time accurate estimation of the unknown parameters. It is a certainty-equivalence adaptive scheme with a least-squares, regulation-based identifier. The scheme was developed for a specific benchmark problem, where the unknown parameters are the reaction coefficient and the high-frequency gain.

Future work may extend the proposed methodology to more complicated cases (systems of parabolic PDEs). Another direction for future research is the development of the certainty-equivalence adaptive scheme with a least-squares, regulation-based identifier to hyperbolic PDEs.

## Appendix

**Proof of Lemma 5.2:** It suffices to notice that the initial-boundary value problem (5.7), (5.8), (5.9) can be written in abstract form as $\dot{w} = Aw + Bw$, where $A$ is a linear operator which is the infinitesimal generator of a $C_0$ semigroup of contractions on $L^2(0,1)$ and $B$ is a linear bounded operator on $L^2(0,1)$. The result follows from Theorem 1.1 on page 76 in [28]. ◁

**Proof of Lemma 5.3:** Let $\gamma \in \Re$ be a constant with $\int_0^1 xk(x)dx + \gamma \|k\|^2 \neq 1$ (notice that for every $k \in C^2([0,1])$ there are infinite values of $\gamma \in \Re$ for which this relation holds) and define

$$\beta := \left(1 - \int_0^1 sk(s)ds - \gamma \|k\|^2\right)^{-1} \tag{A.1}$$

$$a := \beta p \gamma \tag{A.2}$$

Let $\sigma > 0$ be the constant that corresponds to $k \in C^2([0,1])$ being defined by (5.5), (2.6) and to the constant $a, \gamma \in \Re$ defined by (A.2), (A.1) for which (5.10) holds (recall Lemma 5.2). Consider the following Fredholm transformation

$$w(t,x) = u(t,x) - (x + \gamma k(x))\int_0^1 k(s)u(t,s)ds, \text{ for } t \geq 0, \; x \in [0,1] \tag{A.3}$$

and its inverse

$$u(t,x) = w(t,x) + \beta(x + \gamma k(x))\int_0^1 k(s)w(t,s)ds, \text{ for } t \geq 0, \ x \in [0,1]. \quad (A.4)$$

Using (A.1), (A.2), (A.3), (A.4), (5.2), (5.3), (5.4), (5.5), we notice that $w \in C^0(\Re_+; L^2(0,1)) \cap C^1((0,+\infty) \times [0,1])$ with $w[t] \in C^2([0,1])$ for $t > 0$ is a solution of initial-boundary value problem (5.7), (5.8), (5.9) with

$$w_0(x) = u_0(x) - (x + \gamma k(x))\int_0^1 k(s)u_0(s)ds, \text{ for } x \in [0,1] \quad (A.5)$$

Applying the triangular inequality and the Cauchy-Schwartz inequality in (A.4), we get

$$\|u[t]\| \leq \left(1 + \frac{|\beta|}{\sqrt{3}}\|k\| + |\beta\gamma|\|k\|^2\right)\|w[t]\|, \text{ for } t \geq 0 \quad (A.6)$$

Application of the triangular inequality and the Cauchy-Schwartz inequality in (A.5) gives

$$\|w_0\| \leq \left(1 + \frac{1}{\sqrt{3}}\|k\| + |\gamma|\|k\|^2\right)\|u_0\| \quad (A.7)$$

Thus, combining (5.10), (A.6) and (A.7), we obtain estimate (5.6), with

$$Q := \left(1 + \frac{|\beta|}{\sqrt{3}}\|k\| + |\beta\gamma|\|k\|^2\right)\left(1 + \frac{1}{\sqrt{3}}\|k\| + |\gamma|\|k\|^2\right) \quad (A.8)$$

The proof is complete. ◁

**Proof of Lemma 5.4:** Assume that there exists $\bar{u} \in C^0(\Re_+; L^2(0,1)) \cap C^1((0,+\infty) \times [0,1])$ with $\bar{u}[t] \in C^2([0,1])$ for $t > 0$, satisfying (5.1) as well as the boundary conditions (5.2), (5.3) and the initial condition $\bar{u}[0] = u_0 \in L^2(0,1)$. The latter implies that $y = u - \bar{u}$ is also a solution of the initial boundary value problem (5.1), (5.2), (5.3), (5.4) with $y \in C^0(\Re_+; L^2(0,1)) \cap C^1((0,+\infty) \times [0,1])$ with $y[t] \in C^2([0,1])$ for $t > 0$ which satisfies the boundary conditions (5.2), (5.3) as well as the initial condition $y[0] \equiv 0 \in L^2(0,1)$. In addition from Lemma 5.3, it follows that $y$ satisfies (5.6). Hence, we have

$$\|y[t]\| \leq Q\exp(\sigma t)\|y[0]\|, \quad \forall t \geq 0 \quad (A.9)$$

where $Q, \sigma > 0$ are the constants that are guaranteed to exist by virtue of Lemma 5.3. Replacing $y = u - \bar{u}$ in (A.9) and since $y[0] \equiv 0 \in L^2(0,1)$, we get

$$\|u[t] - \bar{u}[t]\| \leq 0, \ \forall t \geq 0 \quad (A.10)$$

Inequality (A.10) implies that $u[t] = \bar{u}[t]$ for all $t \geq 0$. The proof is complete. ◁

**Proof of Lemma 5.5:** Using (5.13), the Cauchy-Schwartz inequality and the fact that $\|\phi_n\| = 1$ for $n = 1, 2, \ldots$, $x\exp(-x) \leq \exp(-1)$ for $x \geq 0$, we get for every $t_0 \in (0, T]$ and $t \in [t_0, T]$:

$$\left|\exp(-pn^2\pi^2 t)c_n\right| \leq \sqrt{2}\|w_0\|\frac{\exp(-1)}{pn^2\pi^2 t_0}, \ n = 1, 2, \ldots \quad (A.11)$$

Furthermore, using (5.14), the fact that $|\exp(At)| \leq \exp(|A|t)$ for all $t \geq 0$ and the Cauchy-Schwartz inequality, we obtain the following estimate for every $t_0 \in (0, T]$ and $t \in [t_0, T]$:

$$\left|L_n \int_0^t \exp(-pn^2\pi^2(t-s))r\exp(As)\xi\, ds\right| \leq |r|\|\xi\|\exp(|A|T)\frac{\sqrt{2}}{pn^3\pi^3}, \ n = 1, 2, \ldots \quad (A.12)$$

Inequalities (A.11) and (A.12) imply that the infinite series appearing on the right hand side of (5.15) are uniformly and absolutely convergent on $[t_0, T] \times [0,1]$, for all $t_0 \in (0, T]$. Therefore, $w \in C^0((0,T] \times [0,1])$. Combining (5.15), (5.16), (5.17), (5.18), the fact that $w_0 \in L^2(0,1)$ and using Parseval's identity we conclude that for any $\varepsilon > 0$, there exists $\delta > 0$ such that

$$\|w[t]-w_0\|\le\varepsilon, \text{ for any } t\in[0,\delta] \tag{A.13}$$

which shows that the mapping $[0,T]\ni t\to w[t]\in L^2(0,1)$ is continuous. Indeed, the fact that $w_0\in L^2(0,1)$ in conjunction with (5.17) implies that $\sum_{n=1}^{\infty}c_n^2<+\infty$. The latter combined with (5.14), (5.15), (5.16), (5.17), (5.18), the fact that $|\exp(At)-I|\le|A|t\exp(|A|t)$ for all $t\ge 0$ and the fact that $\|\phi_n\|=1$ shows that for all $t\in(0,T]$ and for every integer $m\ge N+1$, it holds that

$$\|w[t]-w_0\|^2 \le N|A|^2 t^2 \exp(2|A|T)\sum_{n=1}^{N}c_n^2 + 2\left(1-\exp(-pm^2\pi^2 t)\right)^2 \sum_{n=N+1}^{m}c_n^2$$
$$+2t|r|^2\exp(2|A|T)\left(\sum_{n=1}^{N}c_n^2\right)\sum_{n=N+1}^{\infty}\frac{1}{pn^4\pi^4}+2\sum_{n=m+1}^{\infty}c_n^2 \tag{A.14}$$

Inequality (A.13) is a direct consequence of inequality (A.14) for $m$ sufficiently large and $\delta>0$ sufficiently small.

We show next that $w[t]\in C^1([0,1])$ for every $t\in(0,T]$ and satisfies the following equation for all $(t,x)\in(0,T]\times[0,1]$:

$$\frac{\partial w}{\partial x}(t,x)=\sum_{n=1}^{N}\phi_n'(x)e_n\exp(At)\xi + \sum_{n=N+1}^{\infty}\phi_n'(x)\exp(-pn^2\pi^2 t)c_n$$
$$+\sum_{n=N+1}^{\infty}L_n\phi_n'(x)\int_0^t\exp(-pn^2\pi^2(t-s))r\exp(As)\xi\,ds \tag{A.15}$$

In order to show this, it suffices to show that the infinite series appearing on the right hand side of (A.15) are uniformly and absolutely convergent on $[t_0,T]\times[0,1]$, for all $t_0\in(0,T]$. First, notice that according to definition (2.6), the following estimate holds

$$\max_{0\le x\le 1}(|\phi_n'(x)|)\le n\pi\sqrt{2}, \text{ for all } n=1,2,... \tag{A.16}$$

Combining (5.13), (A.16), the fact that $\|\phi_n\|=1$ for $n=1,2,...$, and the fact that $x^a\exp(-x)\le a^a\exp(-a)$, for any $x\ge 0$, $a>0$ and applying the Cauchy-Schwartz inequality, we get for all $x\in[0,1]$, $t_0\in(0,T]$ and $t\in[t_0,T]$:

$$\left|\phi_n'(x)\exp\left(-pn^2\pi^2 t\right)c_n\right|\le\|w_0\|\frac{1}{2n^2\pi^2}\left(\frac{3}{2pt_0 e}\right)^{3/2}, \ n=1,2,... \tag{A.17}$$

Furthermore, using (A.16), (A.12) and the Cauchy-Schwartz inequality, we obtain

$$\left|L_n\phi_n'(x)\int_0^t\exp\left(-pn^2\pi^2(t-s)\right)r\exp(As)\xi\,ds\right|\le|r||\xi|\exp(|A|T)\frac{\sqrt{2}}{pn^2\pi^2},$$
$$\text{for } n=1,2,... \text{ and } x\in[0,1] \tag{A.18}$$

Inequalities (A.17), (A.18) indicate that for any $T>0$, the infinite series appearing on the right hand side of (A.15) are uniformly and absolutely convergent on $[t_0,T]\times[0,1]$, for all $t_0\in(0,T]$. Therefore $w[t]\in C^1([0,1])$. The proof is complete. ◁

**Proof of Lemma 2.1:** We prove all implications one by one.

Implication (i): If $Q_{n,1}(\tau_{i+1},\mu_{i+1})=0$ for $n=1,2,...$ then definition (2.24) in conjunction with continuity of $g_n(t,s)$ for $t,s\in[\mu_{i+1},\tau_{i+1}]$ (a consequence of definition (2.17) and the fact that $u\in C^0([\mu_{i+1},\tau_{i+1}];L^2(0,1))$) implies that $g_n(t,s)=0$ for $t,s\in[\mu_{i+1},\tau_{i+1}]$. Continuity of the mapping $\tau\to\int_0^1\sin(n\pi x)u(\tau,x)dx$ for $\tau\in[\mu_{i+1},\tau_{i+1}]$ (a consequence of the fact that $u\in C^0([\mu_{i+1},\tau_{i+1}];L^2(0,1))$) and definition (2.17) implies that $\int_0^1\sin(n\pi x)u(\tau,x)dx=0$ for $\tau\in[\mu_{i+1},\tau_{i+1}]$ and $n=1,2,...$. Since the set of functions $\{\phi_n(x) = \sqrt{2}\sin(n\pi x): n = 1, 2, ...\}$ is an orthonormal basis of $L^2(0,1)$, we obtain that $u[t] = 0$ (in the sense of $L^2(0,1)$) for $t \in [\mu_{i+1}, \tau_{i+1}]$. Definition (2.16) implies that $f_n(t,s) = 0$ for $t, s \in [\mu_{i+1}, \tau_{i+1}]$. Since $g_n(t,s) = f_n(t,s) = 0$ for $t, s \in [\mu_{i+1}, \tau_{i+1}]$ and since $c > 0$, we conclude from (2.15) that $j_n(t,s) = 0$ for $t, s \in [\mu_{i+1}, \tau_{i+1}]$. Definitions (2.25) and (2.26) imply that $Q_{n,2}(\tau_{i+1}, \mu_{i+1}) = Q_{n,3}(\tau_{i+1}, \mu_{i+1}) = 0$ for $n = 1, 2, ...$. Therefore, definition (2.27) implies that $S_i = \Re^2$.

Implication (ii): Suppose that $S_i = \Re^2$. Moreover, suppose that on the contrary there exists $n \in \{1, 2, ...\}$ such that $Q_{n,1}(\tau_{i+1}, \mu_{i+1}) \neq 0$. Then definition (2.27) implies that
$$S_i \subseteq \left\{ (\vartheta_1, \vartheta_2) \in \Re^2 : \vartheta_1 = \frac{H_{n,1}(\tau_{i+1}, \mu_{i+1})}{Q_{n,1}(\tau_{i+1}, \mu_{i+1})} - \vartheta_2 \frac{Q_{n,2}(\tau_{i+1}, \mu_{i+1})}{Q_{n,1}(\tau_{i+1}, \mu_{i+1})} \right\},$$
which contradicts the assumption $S_i = \Re^2$. Therefore, we conclude that $Q_{n,1}(\tau_{i+1}, \mu_{i+1}) = 0$ for $n = 1, 2, ...$. Then, working as in the proof of implication (i), we conclude that $u[t] = 0$ (in the sense of $L^2(0,1)$) for $t \in [\mu_{i+1}, \tau_{i+1}]$.

Implication (iii): Suppose that $S_i \neq \{(\theta, c)\}$ and $S_i \neq \Re^2$. By virtue of implication (i) there exists $n \in \{1, 2, ...\}$ such that $Q_{n,1}(\tau_{i+1}, \mu_{i+1}) \neq 0$. Define the index set $I$ to be the set of all $n \in \{1, 2, ...\}$ with $Q_{n,1}(\tau_{i+1}, \mu_{i+1}) \neq 0$. We first show that $U(t) \equiv 0$.

Definition (2.27) implies that $S_i \subseteq \left\{ (\vartheta_1, \vartheta_2) \in \Re^2 : \vartheta_1 = \frac{H_{n,1}(\tau_{i+1}, \mu_{i+1})}{Q_{n,1}(\tau_{i+1}, \mu_{i+1})} - \vartheta_2 \frac{Q_{n,2}(\tau_{i+1}, \mu_{i+1})}{Q_{n,1}(\tau_{i+1}, \mu_{i+1})}, n \in I \right\}$. Since $S_i \neq \{(\theta, c)\}$, we conclude that there exists a constant $\lambda \in \Re$ such that $\frac{Q_{n,2}(\tau_{i+1}, \mu_{i+1})}{Q_{n,1}(\tau_{i+1}, \mu_{i+1})} = \lambda$ for all $n \in I$.

Since $S_i \neq \{(\theta, c)\}$, definition (2.27) implies that $\left(Q_{n,2}(\tau_{i+1}, \mu_{i+1})\right)^2 = Q_{n,1}(\tau_{i+1}, \mu_{i+1}) Q_{n,3}(\tau_{i+1}, \mu_{i+1})$ for all $n \in I$. By virtue of the definitions (2.24), (2.25) and (2.26) and the fact that the Cauchy-Schwarz inequality holds as an equality only when two functions are linearly dependent, we obtain for all $n \in I$ the existence of a constant $\gamma_n \in \Re$ such that $j_n(t,s) = \gamma_n g_n(t,s)$ for $t, s \in [\mu_{i+1}, \tau_{i+1}]$ (notice that $g_n \neq 0 \in L^2\left((\mu_{i+1}, \tau_{i+1})^2\right)$ since $n \in I$ and $Q_{n,1}(\tau_{i+1}, \mu_{i+1}) \neq 0$). Since $\frac{Q_{n,2}(\tau_{i+1}, \mu_{i+1})}{Q_{n,1}(\tau_{i+1}, \mu_{i+1})} = \lambda$ for all $n \in I$, we obtain from definitions (2.24), (2.25) and the fact that $j_n(t,s) = \gamma_n g_n(t,s)$ for $t, s \in [\mu_{i+1}, \tau_{i+1}]$ that $j_n(t,s) = \lambda g_n(t,s)$ for $t, s \in [\mu_{i+1}, \tau_{i+1}]$ and $n \in I$. Definitions (2.17) allow us to conclude that
$$U(t) = -\frac{\lambda}{(-1)^n pn\pi} \int_0^1 \sin(n\pi x) u(t,x) dx \text{ for all } t \in (\mu_{i+1}, \tau_{i+1}) \text{ and } n \in I.$$

If $\lambda = 0$ then $U(t) \equiv 0$. We next assume that $\lambda \neq 0$ and we show again that $U(t) \equiv 0$. For this, it suffices to show that there exists $m \in \{1, 2, ...\}$ with $m \notin I$. Indeed, if $m \notin I$ then (working exactly as in the proof of implication (i)) we get $\int_0^1 \sin(m\pi x) u(\tau, x) dx = 0$ for $\tau \in [\mu_{i+1}, \tau_{i+1}]$. But then we get from (2.15), (2.16) and (2.17) that $j_m(t,s) = 0$ for $t, s \in [\mu_{i+1}, \tau_{i+1}]$, which implies that $U(t) \equiv 0$.

Since (2.15) for $t, s \in [\mu_{i+1}, \tau_{i+1}]$ implies (2.14) for $\tau \in (\mu_{i+1}, \tau_{i+1})$ and since $U(t) = \frac{\lambda}{(-1)^n pn\pi} \int_0^1 \sin(n\pi x) u(t,x) dx$ for all $t \in (\mu_{i+1}, \tau_{i+1})$ and $n \in I$, we obtain (from continuity of the mapping $\tau \to \int_0^1 \sin(n\pi x) u(\tau, x) dx$ for $\tau \in [\mu_{i+1}, \tau_{i+1}]$) that

$$\int_0^1 \sin(n\pi x)u(t,x)dx = \exp\big((\theta - n^2\pi^2 p + \lambda c)(t - \mu_{i+1})\big)\int_0^1 \sin(n\pi x)u(\mu_{i+1},x)dx, \text{ for } t \in [\mu_{i+1}, \tau_{i+1}] \text{ and } n \in I \quad \text{(A.19)}$$

If $m, n \in I$ then we have $\lambda^{-1} p\pi U(t) = \dfrac{1}{(-1)^n n}\int_0^1 \sin(n\pi x)u(t,x)dx = \dfrac{1}{(-1)^m m}\int_0^1 \sin(m\pi x)u(t,x)dx$ for $t \in (\mu_{i+1}, \tau_{i+1})$. Continuity of the mappings $\tau \to \int_0^1 \sin(n\pi x)u(\tau,x)dx$, $\tau \to \int_0^1 \sin(m\pi x)u(\tau,x)dx$ for $\tau \in [\mu_{i+1}, \tau_{i+1}]$ implies that $\dfrac{1}{(-1)^n n}\int_0^1 \sin(n\pi x)u(t,x)dx = \dfrac{1}{(-1)^m m}\int_0^1 \sin(m\pi x)u(t,x)dx$ for $t \in [\mu_{i+1}, \tau_{i+1}]$. If $\int_0^1 \sin(n\pi x)u(\mu_{i+1},x)dx = 0$ or $\int_0^1 \sin(m\pi x)u(\mu_{i+1},x)dx = 0$ then (A.19) and (2.17) would give that $Q_{n,1}(\tau_{i+1}, \mu_{i+1}) = 0$ or $Q_{m,1}(\tau_{i+1}, \mu_{i+1}) = 0$, which contradicts the fact that $m, n \in I$. Thus $\int_0^1 \sin(n\pi x)u(\mu_{i+1},x)dx \neq 0$ and $\int_0^1 \sin(m\pi x)u(\mu_{i+1},x)dx \neq 0$. Using (A.19) and the fact $\dfrac{1}{(-1)^n n}\int_0^1 \sin(n\pi x)u(t,x)dx = \dfrac{1}{(-1)^m m}\int_0^1 \sin(m\pi x)u(t,x)dx$ for $t \in [\mu_{i+1}, \tau_{i+1}]$ as well as the fact that $\int_0^1 \sin(n\pi x)u(\mu_{i+1},x)dx \neq 0$ and $\int_0^1 \sin(m\pi x)u(\mu_{i+1},x)dx \neq 0$, we obtain the equation $\exp\big(-n^2\pi^2 p(t - \mu_{i+1})\big) = \exp\big(-m^2\pi^2 p(t - \mu_{i+1})\big)$ for $t \in [\mu_{i+1}, \tau_{i+1}]$, which can only hold for $m = n$. Therefore, the index set $I$ is a singleton, which implies the existence of $m \in \{1, 2, ...\}$ with $m \notin I$.

Since $U(t) \equiv 0$, it follows from (2.17) that $j_n(t,s) = 0$ for $t, s \in [\mu_{i+1}, \tau_{i+1}]$ and $n = 1, 2, ...$. Consequently, we obtain from (2.25) and (2.26) that $Q_{n,2}(\tau_{i+1}, \mu_{i+1}) = Q_{n,3}(\tau_{i+1}, \mu_{i+1}) = 0$ for $n = 1, 2, ...$. Therefore, definition (2.27) implies that $S_i = \left\{(\vartheta_1, \vartheta_2) \in \Re^2 : \vartheta_1 = \dfrac{H_{n,1}(\tau_{i+1}, \mu_{i+1})}{Q_{n,1}(\tau_{i+1}, \mu_{i+1})}, n \in I\right\}$. Since $(\theta, c) \in S_i$ it follows that $S_i = \left\{(\theta, \vartheta_2) \in \Re^2 : \vartheta_2 \in \Re\right\}$. The proof is complete. ◁

## References

[1] Anfinsen, H. and O. M. Aamo, "Adaptive Stabilization of n+1 Coupled Linear Hyperbolic Systems with Uncertain Boundary Parameters Using Boundary Sensing", *Systems and Control Letters*, 99, 2017, 72-84.
[2] Anfinsen, H. and O. M. Aamo, "Model Reference Adaptive Control of n+1 Coupled Linear Hyperbolic PDEs", *Systems and Control Letters*, 109, 2017, 1-11.
[3] Anfinsen, H. and O. M. Aamo, "Adaptive Output Feedback Stabilization of Coupled Linear Hyperbolic PDEs with Uncertain Boundary Conditions", *SIAM Journal of Control and Optimization*, 55, 2017, 3928-3946.
[4] Anfinsen, H., M. Diagne, O. M. Aamo, and M. Krstic, "An Adaptive Observer Design for n+1 Coupled Linear Hyperbolic PDEs Based on Swapping", *IEEE Transactions on Automatic Control*, 61, 2016, 3979-3990.
[5] Bekiaris-Liberis, N. and M. Krstic, "Delay-Adaptive Feedback for Linear Feedforward Systems", *Systems and Control Letters*, 59, 2010, 277-283.
[6] Bernard, P. and M. Krstic, "Adaptive Output-Feedback Stabilization of Non-Local Hyperbolic PDEs", *Automatica*, 50, 2014, 2692-2699.
[7] Bresch-Pietri, D. and M. Krstic, "Adaptive Trajectory Tracking Despite Unknown Input Delay and Plant Parameters", *Automatica*, 45, 2009, 2075-2081.

[8] Bresch-Pietri, D. and M. Krstic, "Delay-Adaptive Predictor Feedback for Systems with Unknown Long Actuator Delay", *IEEE TAC*, 55, 2010, 2106-2112.

[9] Bresch-Prietri, D. and M. Krstic, "Delay-Adaptive Control for Nonlinear Systems", *IEEE Transactions on Automatic Control*, 59, 2014, 1203-1218.

[10] Bresch-Pietri, D. and M. Krstic, "Output-Feedback Adaptive Control of a Wave PDE with Boundary Anti-Damping", *Automatica*, 50, 2014, 1407-1415.

[11] Bresch-Pietri, D. and M. Krstic, "Adaptive Compensation of Diffusion-Advection Actuator Dynamics Using Boundary Measurements", *Proceedings of the 54th IEEE Conference on Decision and Control*, 2015.

[12] Christophides, P., *Nonlinear and Robust Control of Partial Differential Equation Systems: Methods and Applications to Transport-Reaction Processes*, Birkhäuser, 2001.

[13] El-Farrah, N. and P. Christophides, "Hybrid Control of Parabolic PDE Systems", *Proceedings of the 41st IEEE Conference on Decision and Control*, 2002, 216-22.

[14] Espiritia, N., A. Girard, N. Marchand and C. Prieur, "Event-Based Boundary Control of a Linear 2x2 Hyperbolic System Via Backstepping Approach", to appear in *IEEE TAC*.

[15] Espiritia, N., A. Girard and C. Prieur, "Event-Based Control of Linear Hyperbolic Systems of Conservation Laws", *Automatica*, 70, 2016, 275-287.

[16] Jiang, Z., B. Cui, W. Wu and B. Zhuang, "Event-Driven Observer-Based Control for Distributed Parameter Systems Using Mobile Sensor and Actuator", *Computers and Mathematics with Applications*, 72, 2016, 2854-2864.

[17] Karafyllis, I. and M. Krstic, "Adaptive Certainty-Equivalence Control With Regulation-Triggered Finite-Time Least-Squares Identification", to appear in *IEEE Transactions on Automatic Control* (see also: arXiv:1609.03016 [math.OC] and arXiv:1609.03017 [math.OC]).

[18] Karafyllis, I. and M. Krstic, "Sampled-Data Boundary Feedback Control of 1-D Parabolic PDEs", *Automatica*, 87, 2018, 226-237.

[19] Krstic, M., "Systematization of Approaches to Adaptive Boundary Control of PDEs", *International Journal of Robust and Nonlinear Control*, 16, 2006, 801-818.

[20] Krstic M., "Adaptive Control of an Unstable Wave PDE", *Proceedings of the 2009 American Control Conference*, 2009.

[21] Krstic, M. and A. Smyshlyaev, "Adaptive Control of PDEs", *Annual Reviews in Control*, 32, 2008, 149-160.

[22] Krstic, M. and A. Smyshlyaev, *Boundary Control of PDEs: A Course on Backstepping Designs*, SIAM, 2008.

[23] Krstic, M. and A. Smyshlyaev, "Adaptive Boundary Control for Unstable Parabolic PDEs - Part I: Lyapunov Design", *IEEE TAC*, 53, 2008, 1575-1591.

[24] Lamare, P.-O., A. Girard and C. Prieur, "Switching Rules for Stabilization of Linear Systems of Conservation Laws", *SIAM Journal of Control and Optimization*, 53, 2015, 1599-1624.

[25] Lasiecka, I. and R. Triggiani, "Stabilization and Structural Assignment of Dirichlet Boundary Feedback Parabolic Equations", *SIAM Journal on Control and Optimization*, 21, 1983, 766-803.

[26] Masashi, W. and S. Hideki, "Event-Triggered Control of Infinite Dimensional Systems", arXiv1801.05521v1 [math.OC].

[27] Mechoud, S., "Adaptive Distributed Parameter and Input Estimator in Linear Parabolic PDEs", *International Journal of Adaptive Control and Signal Processing*, 30, 2016.

[28] Pazy, A., *Semigroups of Linear Operators and Applications to Partial Differential Equations*, New York, Springer, 1983.

[29] Prieur, C., A. Girard and E. Witrant, "Stability of Switched Linear Hyperbolic Systems by Lyapunov Techniques", *IEEE Transactions on Automatic Control*, 59, 2014, 2196-2202.

[30] Selivanov, A. and E. Fridman, "Distributed Event-Triggered Control of Diffusion Semilinear PDEs", *Automatica*, 68, 2016, 344-351.

[31] Smyshlyaev, A., and M. Krstic, "Closed-Form Boundary State Feedbacks for a Class of 1-D Partial Integro-Differential Equations", *IEEE Transactions on Automatic Control*, 49, 2004, 2185-2202.

[32] Smyshlyaev, A. and M. Krstic, "Adaptive Boundary Control for Unstable Parabolic PDEs - Part II: Estimation-Based Designs", *Automatica*, 43, 2007, 1543-1556.

[33] Smyshlyaev, A. and M. Krstic, "Adaptive Boundary Control for Unstable Parabolic PDEs - Part III: Output-Feedback Examples with Swapping Identifiers", *Automatica*, 43, 2007, 1557-1564.

[34] Smyshlyaev, A. and M. Krstic, *Adaptive Control of Parabolic PDEs*, Princeton University Press, 2010.